\documentclass[11pt]{article}

\usepackage {amsmath}
\usepackage {hyperref}
\usepackage {amsfonts}
\usepackage {amsthm}
\usepackage[parfill]{parskip}
\usepackage {amssymb}
\usepackage {fullpage}
\usepackage {mathrsfs}
\usepackage {amscd}
\usepackage[all,cmtip]{xy}

\title{Minimal Idempotents on Solvable Groups}
\author{Tanmay Deshpande}
\date{}
\newcommand{\xto}{\xrightarrow}
\newcommand{\onto}{\twoheadrightarrow}
\newcommand{\hZ}{\hat{\Z}}
\newcommand{\hX}{\hat{X}_*}
\newcommand{\Perv}{\hbox{Perv}}

\newcommand{\av}{{\hbox{av}}}
\newcommand{\tMg}{\tM_{Ug,e}}
\newcommand{\tMga}{\tM_{Ug_1,e}}
\newcommand{\tMgab}{\tM_{Ug_1g_2,e}}
\newcommand{\tMgb}{\tM_{Ug_2,e}}

\newtheorem {thm} {Theorem} [section]
\newtheorem {prop} [thm] {Proposition}
\newtheorem {conj} [thm] {Conjecture}
\newtheorem {lem} [thm] {Lemma}
\newtheorem {cor} [thm] {Corollary}
\theoremstyle{definition}
\newtheorem {defn} [thm] {Definition}
\theoremstyle{remark}
\newtheorem {rk} [thm]  {Remark}
\newtheorem {ex} [thm] {Example}

\newcommand{\bex}{\begin{ex}}
\newcommand{\eex}{\end{ex}}
\renewcommand {\bar} {\overline}
\newcommand{\bpf}{\begin{proof}}
\newcommand{\epf}{\end{proof}}
\newcommand{\rar}[1]{\stackrel{#1}{\longrightarrow}}

\newcommand{\h}{\hbox}
\newcommand{\ind}{\h{ind}}
\newcommand{\indg}{\h{ind}_{G'}^G}
\newcommand{\Resg}{\h{Res}_{G'}^G}

\renewcommand{\L}{\mathcal{L}}
\newcommand{\R}{\mathscr{R}}
\renewcommand{\P}{\mathscr{P}}
\newcommand{\Pn}{{\mathscr{P}_{norm}}}
\newcommand{\N}{\mathscr{N}}

\newcommand{\tK} {\widetilde{K}}

\newcommand{\Q} {\mathbb{Q}}

\newcommand{\F} {\mathbb{F}}
\newcommand{\un} {\mathbf{1}}

\newcommand{\FDG} {\mathcal{F}^G}
\newcommand{\iFDG} {\mathcal{{F}'}^G}
\newcommand{\iFD} {\mathcal{F}'}
\newcommand{\normal} {\triangleleft}
\renewcommand{\t}{\widetilde}
\newcommand{\tU}{\widetilde{U}}

\renewcommand{\k} {\mathtt{k}}

\newcommand{\M} {\mathcal{M}}

\newcommand{\D} {\mathscr{D}}
\newcommand{\C} {\mathcal{C}}

\newcommand{\Hom} {\hbox{Hom}}
\renewcommand{\Vec}{\mbox{Vec}}

\newcommand{\pr}{\mbox{pr}}

\newcommand{\Stab}{\mbox{Stab}}

\newcommand{\Qlcl} {\overline{\mathbb{Q}}_l}
\newcommand{\beq}{\begin{equation}}
\newcommand{\eeq}{\end{equation}}
\newcommand{\bthm}{\begin {thm}}
\newcommand{\ethm}{\end {thm}}
\newcommand{\bprop}{\begin {prop}}
\newcommand{\eprop}{\end {prop}}
\newcommand{\bcor}{\begin {cor}}
\newcommand{\ecor}{\end {cor}}
\newcommand{\blem}{\begin{lem}}
\newcommand{\elem}{\end{lem}}
\newcommand{\bdefn}{\begin{defn}}
\newcommand{\edefn}{\end{defn}}
\newcommand{\brk}{\begin{rk}}
\newcommand{\erk}{\end{rk}}

\newcommand{\G} {\mathbb{G}}

\newcommand{\g} {{\gamma}}
\renewcommand{\l} {{\lambda}}

\newcommand{\id}{\hbox{id}}

\newcommand{\tM}{\widetilde{\M}}

\newcommand{\bit}{\begin{itemize}}
\newcommand{\eit}{\end{itemize}}

\newcommand{\avg}{\av_{G/G'}}
\newcommand{\DG}{\D_G(G)}
\newcommand{\DGp}{\D_{G'}(G')}

\newcommand{\cpu}{\mathfrak{cpu}}
\newcommand{\can}{\mathbb{K}}

\newcommand{\cpuc}{\mathfrak{cpu}^\circ}

\newcommand{\E}{\mathcal{E}}
\newcommand{\sE}{\mathscr{E}}
\newcommand{\Z}{\mathbb{Z}}
\renewcommand{\Vec}{\mbox{Vec}}
\newcommand{\bconj}{\begin{conj}}
\newcommand{\econj}{\end{conj}}
\begin{document}
\maketitle
\begin{abstract}
In this paper, we begin to develop a theory of character sheaves on an affine algebraic group $G$ defined over an algebraically closed field $\k$ of characteristic $p>0$ using the approach developed by Boyarchenko and Drinfeld for unipotent groups. Let $l$ be a prime different from $p$. Following Boyarchenko and Drinfeld (\cite{BD08}), we define the notion of an admissible pair on $G$ and the corresponding idempotent in the $\Qlcl$-linear triangulated braided monoidal category $\D_G(G)$ of conjugation equivariant $\Qlcl$-complexes (under convolution with compact support) and study their properties.  In the spirit of \cite{BD08}, we aim to break up the braided monoidal category $\D_G(G)$  into smaller and more manageable pieces corresponding to these idempotents in $\D_G(G)$. Drinfeld has conjectured that the idempotent in $\D_G(G)$ obtained from an admissible pair is in fact a minimal idempotent and that any minimal idempotent in $\D_G(G)$ can be obtained from some admissible pair on $G$. In this paper, we prove this conjecture in the case when the neutral connected component $G^\circ \subset G$ is a solvable group. For general groups, we prove that this conjecture is in fact equivalent to an apriori weaker conjecture. Using these results, we reduce the problem of defining character sheaves on general algebraic groups to a special case which we call the ``Heisenberg case''. Moreover as we will see in this paper, the study of character sheaves in the Heisenberg case may be considered, in a certain sense, as a twisted version of the theory of character sheaves on reductive groups as developed by Lusztig in \cite{L}.
\end{abstract}

\tableofcontents

\section{Introduction}
Let $G$ be an affine algebraic group over an algebraically closed field $\k$ of characteristic $p>0$ and let $l$ be a prime number not equal to $p$. The goal of the theory of character sheaves is to geometrize the notion of irreducible characters of finite groups to the setting of algebraic groups.  One of the motivations to the theory of character sheaves comes from the case when $\k$ is an algebraic closure of a finite field $\F_q$ and the group $G$ has an $\F_q$-structure. In this situation, we would like to study the irreducible characters of the finite groups $G(\F_{q^n})$ for positive integers $n$ in terms of the geometry of $G$. For this, it is convenient to fix a prime $l\neq p$ and study $\Qlcl$-valued irreducible characters of the finite groups $G(\F_{q^n})$ since this allows us to use the theory of $\Qlcl$-sheaves and complexes and Grothendieck's sheaf-function correspondence.  In the series of papers \cite{L}, Lusztig developed a theory of character sheaves on reductive groups and related it to the character theory of finite groups of Lie type. Later, Lusztig conjectured the existence of an interesting theory of character sheaves on unipotent groups. Inspired by Lusztig's work, Boyarchenko and Drinfeld (in \cite{BD06},\cite{BD08}, \cite{Bo13})  developed a theory of character sheaves on unipotent groups and related it to the character theory of finite unipotent groups. We refer to \cite{BD06} for a more detailed exposition about the motivation behind the theory of character sheaves on algebraic groups. This paper is a first step towards developing a theory of character sheaves on general algebraic groups. 

\subsection{Notation and conventions}
Throughout this paper, we fix the two distinct primes $p, l$ as well as the algebraically closed field $\k$ of characteristic $p>0$.  All algebraic groups and schemes will be assumed to be over the field $\k$, unless explicitly stated otherwise. We will often deal with certain $\Qlcl$-linear triangulated categories and functors associated with schemes and morphisms of schemes. All such functors should always be considered in the derived sense even though we never decorate them with ``$L$'' or ``$R$''. For a $\k$-scheme $X$,  let $\D(X):=D^b_c(X,\Qlcl)$, the $\Qlcl$-linear triangulated {\it symmetric monoidal} category of $\Qlcl$-complexes on $X$ under (derived) tensor product. The constant rank 1 local system $(\Qlcl)_X\in \D(X)$ is the unit object for the symmetric monoidal structure defined by the tensor product. For any integer $n$, we have the degree shift by $n$ functor $(\cdot)[n]:\D(X)\to\D(X)$ as well as the $n$-th Tate twist functor $(\cdot)(n):\D(X)\to\D(X)$. We think of $\D(X)$ as the geometric analogue of the algebra of functions on a finite set under pointwise product. If an algebraic group $G$ acts on $X$, then $\D_G(X)$ denotes the $\Qlcl$-linear category of $G$-equivariant complexes. (See \S\ref{eqdercat} below.) We have the degree shift and Tate twist functors on $\D_G(X)$ as well. 

By an algebraic group, we mean a smooth affine group scheme over $\k$. For an algebraic group $G$, the $\Qlcl$-linear triangulated category $\D(G)$ is also monoidal under convolution with compact support. We recall (from \cite[\S 1.4]{BD08}) that for $M,N\in\D(G)$, we define the convolution with compact support as $M\ast N:=\mu_!(M\boxtimes N)$ where $\mu:G\times G\to G$ is the group multiplication and $\mu_!$ denotes (derived) pushforward with compact support. Here the external tensor product $M\boxtimes N:=p_1^*M\otimes p_2^*N$ where $p_i:G\times G\to G$ for $i=1,2$ are the two projections. The unit object in $\D(G)$ is $\delta_1$, the delta-sheaf supported at $1\in G$ with stalk $\Qlcl$. We also have the $\Qlcl$-linear triangulated {\it braided monoidal category} $\D_G(G)$ of $G$-equivariant complexes for the conjugation action of $G$ on itself with unit object $\delta_1$.  For $M\in\D(G)$ and $N\in \D_G(G)$, we have braiding isomorphisms $\beta_{M,N}:M\ast N\to N\ast M$ defined in \cite[Defn. A.43]{BD08} which provide the braided structure on $\D_G(G)$. The category $\D_G(G)$ is equipped with a twist $\theta$, an automorphism of the identity functor on $\D_G(G)$. This gives $\D_G(G)$ the structure of a ribbon $\mathfrak{r}$-category. We refer to \cite[\S5.5]{BD06}, \cite[Appendix B]{BD08} for more on the structure of the categories $\D(G)$ and $\D_G(G)$. The monoidal category $\D(G)$ can be thought of as the geometric analogue of the algebra of functions on a finite group under convolution, or equivalently the group algebra of a finite group. The braided monoidal category $\D_G(G)$ is the geometric analogue of the (commutative)  algebra of class functions on a finite group under convolution.

\subsection{Main goal}
Let $G$ be any affine algebraic group over $\k$. Character sheaves on $G$ are supposed to be some special objects in the category $\DG$. Informally, our goal is to break the braided triangulated category $\D_G(G)$ into ``small and manageable pieces'' corresponding to certain idempotents in $\D_G(G)$ and to then define character sheaves in each such piece. For motivation, let us consider the case studied by Boyarchenko and Drinfeld where $G$ is unipotent. Let $e\in \D_G(G)$ be a minimal idempotent (see Definition \ref{defmi}). Let us consider the Hecke subcategory $e\DG\subset \DG$. It is proved in \cite{BD08} that with the minimal idempotent $e$, we can associate a modular category $\M_{G,e}\subset e\D_G(G)$ such that $e\D_G(G)\cong D^b(\M_{G,e})$. The simple objects of the modular category $\M_{G,e}\subset e\DG$ (up to a suitable degree shift)  are defined to be the character sheaves on $G$ associated with the minimal idempotent $e\in \DG$.  Also in {\it op. cit.}, the notion of an admissible pair for a unipotent group is defined and it is proved that every minimal idempotent in $\D_G(G)$ can be obtained from a (possibly non-unique) admissible pair for $G$ using a certain induction with compact support functor. Following a suggestion by Drinfeld, we define the notion of admissible pairs on not necessarily unipotent algebraic groups. Many of the results and proofs related to admissible pairs and induction functors from {\it op. cit.} easily generalize to arbitrary algebraic groups after minor modifications.  Drinfeld has conjectured that with every admissible pair for an algebraic group $G$, we can associate a minimal idempotent in $\D_G(G)$ and that every minimal idempotent in $\D_G(G)$ can be obtained in this way. In this paper, we prove this conjecture in the case of groups $G$ whose neutral connected component $G^\circ$ is solvable. We also reduce the problem of defining and studying character sheaves on general algebraic groups to a special case known as the ``Heisenberg case''. We hope that this approach will eventually lead to an interesting theory of character sheaves on arbitrary affine algebraic groups. 

\subsection*{Acknowledgments}
I would like to thank Vladimir Drinfeld for introducing me to the theory of character sheaves on unipotent and reductive groups. Many of the results and conjectures in this paper were conjectured by Drinfeld. Many of the ideas in this paper are based on the works of Boyarchenko and Drinfeld on unipotent groups. I would like to thank Masoud Kamgarpour for many useful discussions which were possible thanks to the support from the University of Queensland, Australia. I would like to thank Antonio Rojas Le\'on for useful correspondence. This work was supported by World Premier Institute Research Center Initiative (WPI), MEXT, Japan.

\section{Definitions, main results and conjectures}\label{summary}
In this section we define some of the main notions discussed in this paper and state the main results proved in this paper and some conjectures.

\subsection{Perfect groups and schemes}
Although we want to study affine algebraic groups $G$, it is often more convenient for us to work in the setting of perfect quasi-algebraic groups (i.e. perfectizations of affine algebraic groups) and perfect quasi-algebraic schemes (i.e. perfectizations of schemes of finite type) over $\k$. (See \cite[\S1.9]{BD08} for more details.) This is because we use the notion of Serre dual (see \S\ref{serredual}) of a connected unipotent group which is only defined canonically as a perfect group scheme. Hence from now on, unless stated otherwise, we will assume that $G$ is the perfectization of an affine algebraic group and that all schemes we consider are perfectizations of schemes of finite type over $\k$. For example, by an abuse of notation, we will continue to use $\G_a$ to denote the perfectization of the additive group.  Even though we work in the setting of perfect quasi-algebraic groups in this paper, most of our constructions and results are applicable for usual affine algebraic groups as well.

We remark that the \'{e}tale topos of a scheme remains unchanged after passing to its perfectization and we can carry out all constructions about derived categories of constructible complexes in the same way for perfect schemes.

\subsection{Equivariant derived categories}\label{eqdercat}
Let $G$ be a perfect quasi-algebraic group acting on a perfect quasi-algebraic scheme $X$. Let $\alpha:G\times X\to X$ be the action map and let $p:G\times X\to X$ denote the second projection. 

We define the equivariant derived category by setting $\D_G(X):=\D(G\backslash X)$, the $\Qlcl$-linear triangulated category of $\Qlcl$-complexes on the quotient stack $G\backslash X.$ Note that we have the forgetful functor $\D_G(X)\to \D(X)$ and by an abuse of notation, we will often use the same symbol to denote an object of $\D_G(X)$ and its image after applying the forgetful functor.

Let $X,Y$ be perfect quasi-algebraic schemes equipped with actions of $G$ and let $f:X\to Y$ be a $G$-morphism. 
The functors $f^*,f^!,f_*,f_!$ lift canonically to functors between $\D_G(X)$ and $\D_G(Y)$. Verdier duality and tensor product also lift canonically to the equivariant derived categories.

Sometimes it is useful to think of a certain na\"ive equivariant category. Let us define the category $\D_G^{naive}(X)$ whose objects are pairs $(M,\phi)$ where $M$ is an object of $\D(X)$ and $\phi:\alpha^*M\xto\cong p^*M$ is an isomorphism satisfying certain cocycle condition (see \cite[Defn. 1.3]{BD08} for details). The isomorphism $\phi$ is known as a $G$-equivariant structure. We have a forgetful functor $\D_G^{naive}(X)\to \D(X)$. In general $\D_G^{naive}(X)$ may not be a triangulated category. The forgetful functor $\D_G(X)\to \D(X)$ factors as 
\beq
\D_G(X)\to\D_G^{naive}(X)\to\D(X).
\eeq
In particular each object of $\D_G(X)$ has a natural $G$-equivariant structure. In case $G$ is such that $G^\circ$ is unipotent,  the functor $\D_G(X)\to\D_G^{naive}(X)$ is an equivalence.

The triangulated category $\D_G(X)$ has a perverse $t$-structure whose heart is denoted by $\Perv_G(X)$. We have $\Perv_G(X)\cong\Perv_G^{naive}(X)$, where the latter category is defined as the full subcategory of $\D_G^{naive}(X)$ consisting of pairs $(M,\phi)\in \D_G^{naive}(X)$ with $M$ perverse.

If $G$ acts on a perfect quasi-algebraic group $H$ by group automorphisms, then $\D_G(H)$ is a monoidal category under convolution (with compact support). If $G$ acts on groups $H_1,H_2$ by group automorphisms and if $\phi:H_1\to H_2$ is a $G$-morphism, then $\phi_!:\D_G(H_1)\to \D_G(H_2)$ is a monoidal functor.

If we let $G$ act on itself by conjugation, $\D_G(G)$ has the structure of a braided monoidal category. (See \cite[Appendix B]{BD08}.) In fact, for each $M\in \D_G(G)$ and $N\in \D(G)$ we have braiding isomorphisms $\beta_{M,N}:M\ast N\to N\ast M$ (after applying the forgetful functor to $M$).

\subsection{Idempotents in monoidal categories}
Since idempotents in the braided monoidal category $\D_G(G)$ play a key role in our approach to character sheaves, let us recall some notions related to idempotents in monoidal categories. We refer to \cite[\S2]{BD08} and \cite{BDx} for details.

\subsubsection{Open and closed idempotents}
\bdefn
Let $(\M,\otimes,\mathbf{1})$ be a monoidal category.
\bit
\item[(i)] An object $e\in \M$ is said to be an idempotent if $e\otimes e\cong e$.
 
\item[(ii)] A morphism $\un {\xrightarrow{\pi}} e$ is said to be a closed idempotent arrow if both
$$\un\otimes e\xto{\pi\otimes \id_e}e\otimes e \mbox{   and   } e\otimes \un\xto{\id_e\otimes \pi}e\otimes e$$
are isomorphisms. An object $e\in \M$ is said to be a closed idempotent if there exists a closed idempotent arrow $\un\to e$.
\item[(iii)] A morphism $e {\xrightarrow{\pi}} \un$ is said to be an open idempotent arrow if both
$$e\otimes e\xto{\pi\otimes \id_e}\un\otimes e \mbox{   and   } e\otimes e\xto{\id_e\otimes \pi}e\otimes \un$$
are isomorphisms. An object $e\in \M$ is said to be an open idempotent if there exists an open idempotent arrow $e\to \un$.

\eit
\edefn

If $e\in\M$ is an idempotent, the full subcategory $e\M\subset \M$ is defined as the essential image of the functor $\M\to \M$ defined by $X\mapsto e\otimes X.$ We have $e\M=\{X\in \M|e\otimes X\cong X\}$. Similarly we can define $\M e$ and $e\M e$. If $e$ is either an open or a closed idempotent, then in fact $e\M e$ is a monoidal category with unit $e$. In general $e\M e$ is only a semigroupal category. The semigroupal category $e\M e$ is known as the Hecke subcategory associated with the idempotent $e$.

\blem\label{cloinclo} (See \cite[Lem. 1.11]{BDx}.)
Let $\un\to e$ (resp. $e\to \un$) be a closed (resp. open) idempotent arrow in $\M$ and let $e\to f$ (resp. $f\to e$) be any arrow in $\M$. Then  $e\to f$ (resp. $f\to e$) is a closed (resp. open) idempotent in the monoidal category $e\M e$ if and only if the composition, $\un\to f$ (resp. $f\to \un$) is a closed (resp. open) idempotent in $\M$.
\elem

\subsubsection{Locally closed idempotents in triangulated monoidal categories}

Our next goal is to define locally closed idempotents. We have two ways of defining locally closed idempotents in general. We may either define them as open idempotents inside the Hecke subcategory associated with a closed idempotent, or as closed idempotents inside the Hecke subcategory associated with an open idempotent. However as we will see below, both these are equivalent in the case of triangulated monoidal categories.

\bdefn
A triangulated monoidal category is a monoidal category $(\M,\otimes,\un)$ such that $\M$ is a triangulated category and  $\otimes$ is a bitriangulated functor.
\edefn

\blem\label{trsubcat}
If $e$ is either an open or a closed idempotent in a triangulated monoidal category $\M$, then $e\M, \M e, e\M e$ are triangulated subcategories of $\M.$
\elem

\blem\label{opidclid}
Let $e\xto{\pi} \un \xto{\pi'} e'\to e[1]$ be a distinguished triangle in a triangulated monoidal category. Then $\pi$ is an open idempotent if and only if $\pi'$ is a closed idempotent; in which case we have $e\otimes e'=0=e'\otimes e$. This sets up a bijection between the sets of isomorphism classes of closed idempotents and open idempotents in a triangulated monoidal category.
\elem
\bpf
Suppose $\pi$ is an open idempotent. Tensoring the distinguished triangle by $e$ on both sides, we see that $e'\otimes e=0=e\otimes e'.$  Now tensoring the original distinguished triangle by $e'$ on both sides shows that $\pi'$ is a closed idempotent arrow. The other implication can be proved in the same way.
\epf

\bthm\label{lcmain}
Let $\M$ be a triangulated monoidal category and let $f\in \M$ be any object. Then the following are equivalent:
\bit
\item[(i)] There is a closed idempotent $\un \to e$ in $\M$ and an open idempotent $f\to e$ in the Hecke subcategory $e\M e$.
\item[(i$'$)] There is a diagram 
$$\xymatrix{
& \un\ar[d]\\
f\ar[r] & e\\
}$$
where the vertical arrow is a closed idempotent and  both arrows become isomorphisms after tensoring the diagram by $f$ on either side.
\item[(i$''$)] There is a commutative diagram 
$$\xymatrix{
& \un\ar[d]_{\pi}\ar[rd]^{\pi'}& & \\
f\ar[r] & e\ar[r] & e'\ar[r] & f[1]\\
}$$
where $\pi,\pi'$ are closed idempotents in $\M$ and  where the bottom row is a distinguished triangle.
\item[(ii)] There is an open idempotent $u\to \un$ in $\M$ and a closed idempotent $u\to f$ in the Hecke subcategory $u\M u$.
\item[(ii$'$)] There is a diagram 
$$\xymatrix{
u\ar[r]\ar[d] & \un\\
f & \\
}$$
where the horizontal arrow is an open idempotent and  both arrows become isomorphisms after tensoring the diagram by $f$ on either side.
\item[(ii$''$)] There is a commutative diagram 
$$\xymatrix{
u'\ar[d]\ar[rd]^{\psi'}&\\
u\ar[r]^{\psi}\ar[d] & \un\\
f\ar[d] & \\
u'[1]&\\
}$$
where $\psi,\psi'$ are open idempotents and  the first column is a distinguised triangle.
\eit
If $f$ satisfies these equivalent conditions, then $f$ is an idempotent in the triangulated monoidal category $\M$ and the Hecke subcategory $f\M f$ is a triangulated monoidal with unit object $f$.
\ethm
\bpf
The equivalences (i)$\iff$(i$'$)$\iff$(i$''$) and (ii)$\iff$(ii$'$)$\iff$(ii$''$) are straightforward to check using Lemmas \ref{cloinclo} and \ref{opidclid}. Let us now prove that (i$'$) $\implies$ (ii$'$). Given (i$'$), form the diagram 
$$\xymatrix{
u\ar@{-->}[d]\ar[r] & \un\ar[d]_{\pi}\ar[r]^{\pi'}& e'\ar@{=}[d]\ar[r]& u[1]\ar@{-->}[d]\\
f\ar[r] & e\ar[r] & e'\ar[r] & f[1]\\
}$$
using the axioms of a triangulated category by completing the distinguished triangles in the top and bottom row. Now since $f\to e$ is an open idempotent in $e\M e$ and the bottom row is a distinguished triangle, it follows that $e'\in e\M e$ and that $e\to e'$ is a closed idempotent and $f\otimes e'=0=e'\otimes f$ by Lemma \ref{opidclid}. Hence $\un\xto{\pi'} e'$ is a closed idempotent in $\M$ and $u\to \un$ is an open idempotent. Now tensoring the diagram above with $f$ on both sides and using (i$'$) we deduce (ii$''$). The proof that (ii$'$) $\implies$ (i$'$) is similar. The final statement follows by using properties (i), (i$'$) and Lemma \ref{trsubcat}.
\epf

\bdefn\label{lcid}
An object $f\in \M$ satisfying the equivalent conditions in Theorem \ref{lcmain} is said to be a locally closed idempotent in the triangulated monoidal category $\M.$ The Hecke subcategory associated with such an $f$ is the triangulated monoidal category $f\M f$.
\edefn

Using Lemma \ref{cloinclo} and Theorem \ref{lcmain} we get
\bcor\label{lcinlc}
If $e$ is a locally closed idempotent in the triangulated monoidal category $\M$ and if $f$ is a locally closed idempotent in the triangulated monoidal Hecke subcategory $e\M e$, then $f$ is a locally closed idempotent in $\M$.
\ecor

The following remark justifies the use of the adjectives open, closed and locally closed for idempotents:
\brk
Let $j:Y\hookrightarrow X$ be a locally closed (resp. open, closed) embedding of perfect quasi-algebraic schemes. Then it is easy to check that $j_!(\Qlcl)_Y$ is a locally closed (resp. open, closed) idempotent in $(\D(X),\otimes)$, where $(\Qlcl)_Y$ is the constant local system on $Y$.
\erk

\subsubsection{Minimal idempotents in braided monoidal categories}
If $\M$ is a braided monoidal category and if $e\in \M$ is an idempotent, then all the subcategories $e\M, \M e, e\M e$ coincide, and we usually refer to this Hecke subcategory just by $e\M.$ If $e,e'$ are idempotents in a braided monoidal category then so is $e\otimes e'$.

\bdefn\label{defmi}
Let $\M$ be a braided monoidal category with a zero object. An idempotent $e\in \M$ is said to be a minimal idempotent if $e\neq 0$ and if $e'\in \M$ is any idempotent, then either $e\otimes e'\cong e$ or $e\otimes e'=0$.  Equivalently, an idempotent $e\in \M$ is minimal if and only if $e\neq 0$ and the semigroupal Hecke subcategory $e\M$  has no idempotents apart from $0$ and $e$ up to isomorphism.
\edefn

\subsection{Serre duality for connected unipotent groups}\label{serredual}
Let $H$ be a {\it connected} perfect quasi-algebraic group. An object $\L\in\D(H)$ is said to be multiplicative if $\mu^*\L\cong \L\boxtimes \L$, where  $\mu:H\times H\to H$ is the multiplication. It follows that such an $\L$ has to be a rank 1 $\Qlcl$-local system on $H$. Hence we will call such multiplicative objects as {\it multiplicative local systems} on $H$.

In this section we study families of multiplicative local systems on a perfect connected unipotent group $H$. In the commutative case this was studied in \cite{Be} and was later extended to the non-commutative case in \cite[Appendix A]{B}, \cite[\S3]{BD08}. Let us recall here some definitions and facts from {\it loc. cit.} related to this notion. For (connected) unipotent $H$, the ``moduli space'' of multiplicative local systems on $H$ is well defined as a perfect commutative (possibly disconnected) unipotent group $H^*$ (which is known as the Serre dual of $H$) and an $H^*$-family $\sE$ (a universal family) of multiplicative local systems on $H$. The commutative group structure on $H^*$ comes from the fact that the tensor product of multiplicative local systems is again a multiplicative local system. The universal local system $\sE$ is a rank one local system on $H\times H^*$ whose restriction to $H\times \{\L\}\subset H\times H^*$ can be identified with the multiplicative local system $\L$ on $H$, where we consider the multiplicative local system $\L$ as a closed point of $H^*$. If $S$ is a perfect quasi-algebraic variety, an $S$-family of multiplicative local systems on $H$ corresponds to a morphism $S\to H^*$. (See \cite[\S3.1]{BD08} for details.)

Let $c:H\times H\to H$ denote the commutator map. For a multiplicative local system $\L$ on $H$, $c^*\L$ is the constant sheaf $\Qlcl$. However, $\L|_{[H,H]}$ may not be trivial, i.e. the restriction homomorphism $H^*\to [H,H]^*$ is not necessarily the trivial map. However, it is true that the image of this map is finite. (See \cite{Kam} for more.) 
\bprop\label{serredualdecomp}
Let $\E\in [H,H]^*$ lie in the (finite) image of the restriction homomorphism $H^*\to [H,H]^*$. Define $H^*_\E\subset H^*$ to be the fiber of the restriction homomorphism over $\E$, namely the collection of multiplicative local systems on $H$ whose restriction to the commutator subgroup is isomorphic to $\E$. Then
\bit
\item[(i)] We have the pullback isomorphism $(H^{ab})^*\stackrel{\cong}{\longrightarrow} H^*_{\Qlcl}\subset H^*$ induced by the quotient map $$H\onto H^{ab}=H/[H,H].$$
\item[(ii)] $H^{*\circ}=H^*_{\Qlcl}$.
\item[(iii)] $H^*=\coprod\limits_{\E}{H^*_\E}$ is the decomposition of $H^*$ into connected components.
\eit
\eprop

\subsection{Closed idempotents associated with multiplicative local systems}
In this section we study certain special closed idempotents in $\D(G)$ defined using multiplicative local systems on {\it connected unipotent} subgroups of $G$. We recall some notation from \cite[\S4]{BD08}. Let $\P(G)$ denote the set of pairs $(H,\L)$ where $H$ is a connected unipotent subgroup of $G$ and $\L$ a multiplicative local system on $H$. Let $\Pn(G) \subset \P(G)$ denote the subset of pairs $(H,\L)$ such that $H$ is also normal in $G$. For $(H,\L)\in \P(G)$ we can define the normalizer $\N(H,\L)=\N_G(H,\L)$ of the pair (see {\it op. cit.}). 
\bdefn
We define a partial ordering on $\P(G)$ and $\Pn(G)$ by setting $(H_1,\L_1)\leq (H_2,\L_2)$ if $H_1\subset H_2$ and $\L_2|_{H_1}\cong \L_1$. 
\edefn
\bdefn
We say that a pair $(H,\L)$ is compatible with an object $M\in\D(G)$ if $\L\ast M\neq 0$, where $\L$ is considered as an object of $\D(G)$ by extending by zero. We say that pairs $(H_1,\L_1), (H_2,\L_2)\in \P(G)$ are compatible if $\L_1\ast \L_2\neq 0$.  
\edefn

The following lemma is easy to check:
\blem
Pairs $(H_1,\L_1), (H_2,\L_2)\in \P(G)$ are compatible if and only if $$\L_1|_{(H_1\cap H_2)^\circ}\cong \L_2|_{(H_1\cap H_2)^\circ}.$$
\elem

We can define certain closed idempotents using pairs:
\blem\label{cloidepair} (See \cite[Prop. 8.1(a)]{B}.)
Let $(H,\L)\in \P(G)$ be a pair with normalizer $G'$. Then $e_\L=e_{H,\L}:=\L\otimes \can_H\in \D_{G'}(H)$ is a closed idempotent when considered as an object of any of the monoidal categories $\D_{G'}(H)\subset \D_{G'}(G')\subset \D_{G'}(G)$. Here $\can_H$ is the dualizing complex on $H$. It has a canonical $G'$-equivariant structure.
\elem

In our situation the dualizing sheaf $\can_H$ is isomorphic to $\Qlcl[2\dim H]$.

If $(H,\L)\in \P(G)$, we can consider the category $\D_{H,\L}(G)$ of quasi-equivariant complexes (see \cite[\S3.2]{De}) on $G$ with respect to the multiplication by $H$ on the left. Note that for any $M\in\D(G), h\in H$ and $g\in G$,  we have the stalk $(e_\L\ast M)_{hg}\stackrel{\cong}{\longrightarrow} \L_h\otimes (e_\L\ast M)_{g}$ and hence we have an equivalence of categories
\beq
e_\L\D(G)\cong \D_{H,\L}(G).
\eeq



\subsection{Admissible pairs}

Following Drinfeld, let us define the notion of an {\it admissible pair} for $G$. This is similar to the notion of admissible pairs for unipotent groups defined in $\cite{BD08}$.
\bdefn\label{defadm}
An {\it admissible pair} for $G$ is a pair $(H,\L)$ consisting of a connected unipotent subgroup $H\subset G$ and a multiplicative local system $\L$ on $H$ such that the following conditions are satisfied:
\bit
\item[(i)] Let $G'=\N_G(H,\L)$ be the normalizer of $(H,\L)$ in $G$ and let $\R_u(G')$ denote its unipotent radical. Then $\R_u(G')/H$ is commutative.
\item[(ii)] The homomorphism $\phi_\L:\R_u(G')/H\to (\R_u(G')/H)^*$ obtained in this situation (see Remark \ref{auxcon} below) is an isogeny.
\item[(iii)] (Mackey condition) For $g\in G-G'$ we have $$\L|_{(H\cap { }^{g}H)^{\circ}} \ncong { }^{g}\L|_{(H\cap { }^{g}H)^{\circ}}.$$
\eit
In this situation, we have an induced action of $G'$ on the commutative group $\R_u(G')/H$. We say that the admissible pair is {\it central} if this action is trivial.
\edefn

\brk\label{auxcon} (See also \cite[Appendix A.13]{B}.)
If $\R_u(G')/H$ is commutative, we have the commutator map $$c:\R_u(G')\times \R_u(G')\to H.$$ By pulling back the multiplicative local system $\L$ on $H$ along the commutator map we obtain a rank 1 local system on $\R_u(G')\times \R_u(G')$. This can be used to define the homomorphism $\phi_\L:\R_u(G')/H\to (\R_u(G')/H)^*$. Since $\L$ is $G'$-equivariant, $\phi_\L$ is also a $G'$-equivariant map. It in fact defines a $G'$-equivariant skew-symmetric biextension. We refer to \cite[Appendix A]{B} for a detailed exposition of these notions.
\erk

\bdefn
An admissible pair $(H,\L)$ is said to be a {\it Heisenberg} admissible pair if $\N_G(H,\L)=G$. 
\edefn

\brk
If $\N_G(H,\L)=G$, the Mackey condition (iii) is vacuous. Also any admissible pair $(H,\L)$ for $G$ is Heisenberg admissible for $G'=\N_G(H,\L).$ 
\erk

\subsection{Idempotents defined using admissible pairs}
Let $(H,\L)$ be an admissible pair for $G$ with normalizer $G'\subset G$. We have the closed idempotent $e_{\L}=\L\otimes \can_H\in \D_{G'}(G')$. We will define the induction with compact support functor $$\indg:\DGp\to \DG$$ in \S\ref{indfunc}. 
\bdefn\label{idadmpair}
With the above notation, define $f_\L:=\ind_{G'}^G(e_{\L})\in \D_G(G)$. It follows from Proposition \ref{weakid} that $f_\L$ is an idempotent in $\D_G(G)$. This is the idempotent in $\D_G(G)$ associated with the admissible pair $(H,\L)$.
\edefn
\bdefn
If $(H,\L)$ is a Heisenberg admissible pair i.e. if $G'=G$, then $f_\L=e_\L$ is known as a Heisenberg idempotent in $\D_G(G)$. It is clear that a Heisenberg idempotent in $\D_G(G)$ is closed.
\edefn
These idempotents play a key role in our approach to the theory of character sheaves on $G$. One of our goals in this paper is to study these idempotents $f_\L$ and their associated Hecke subcategories $f_\L\DG$.

\subsection{Main conjectures}
Let us now state the main conjectures related to the idempotents $f_\L\in \DG$ associated with admissible pairs $(H,\L)$ that were defined in the previous section. These conjectures are due to Drinfeld following his work on character sheaves for unipotent groups.

\bconj\label{conjmain}
\bit
\item[(i)] The idempotent $f_\L\in \D_G(G)$ defined as above using an admissible pair $(H,\L)$ is locally closed and we have a braided equivalence of braided monoidal categories 
\beq
\indg: e_\L\D_{G'}(G'){\xrightarrow{\cong}} f_\L\DG, 
\eeq
using the same notation as before.
\item[(ii)] The idempotent $f_\L$ associated with an admissible pair $(H,\L)$ is a minimal idempotent in $\D_G(G)$.
\item[(iii)] If $f\in \D_G(G)$ is a minimal idempotent, then $f\cong f_\L$ for some admissible pair $(H,\L)$ for $G$.
\eit
\econj

Eventually our goal is to study the Hecke subcategory $f\DG\subset \DG$ for a minimal idempotent $f\in \DG$ and to define a set $CS_f(G)$ of certain special isomorphism classes of objects in $f\DG$. We will call $CS_f(G)$ the set of character sheaves on $G$ associated with the (minimal) idempotent $f$. If we can prove the above conjectures, then by (iii) and (i), in order to study $f\DG$, it suffices to study the braided monoidal category $e_\L\DGp$, or in other words, it will be sufficient to study Heisenberg idempotents.

\brk
In \cite{BD08}, Boyarchenko and Drinfeld prove these conjectures for {\it unipotent groups} and develop the theory of character sheaves on unipotent groups. In fact it is this work on unipotent groups that led Drinfeld to formulate the above conjectures for general algebraic groups. 
\erk

We now state a conjecture which is apriori weaker than Conjecture \ref{conjmain}.
\bconj\label{conjweak}
Let $(H,\L)$ be a Heisenberg admissible pair for an algebraic group $G$. Then the Heisenberg idempotent $e_\L\in \DG$ is a minimal idempotent.
\econj

\brk
This conjecture is a special case of Conjecture \ref{conjmain}(ii). Moreover, it is not hard to see that this conjecture is in fact equivalent to Conjecture \ref{conjmain}(ii) (see Proposition \ref{weakid}(vii)). We will prove that in fact Conjecture \ref{conjweak} implies all of Conjecture \ref{conjmain} (see Theorem \ref{main4}).
\erk

\subsection{Main results}
One of the main goals of this paper is to prove Conjecture \ref{conjmain} in case $G$ is such that a $G^\circ$ is a {\it solvable group.}

First we develop some general theory for all affine (perfect quasi-) algebraic groups $G$. In \S\ref{avgindfun}, we recall the notion of averaging and induction functors and prove some properties. We also recall the notion of equivariant Fourier-Deligne transform and its compatibility with averaging. In \S\ref{idandind} we construct some idempotents in $\D_G(G)$ using the induction and averaging functors and study their properties. We study the properties of the induction functor in more detail. In \S\ref{anadmpair}, we study admissible pairs in detail. If $(H,\L)$ is an admissible pair for $G$, we classify all admissible pairs $(L,\L')\geq (H,\L)$ and prove that in this case the associated idempotents $f_\L$ and $f_{\L'}\in \DG$ are isomorphic. We also show that if $N\in\D(G)$ is non-zero, then there exists an admissible pair $(H,\L)$ compatible with $N$.

In \S\ref{anheisid} we study Heisenberg idempotents and their associated Hecke subcategories. Let $(H,\L)$ be a Heisenberg admissible pair for $G$ and let $U$ be the unipotent radical of $G$. Then it is clear that $e_\L$ is a Heisenberg idempotent on the connected unipotent group $U$. In \cite{De}, it is shown that $e_\L\D_U(U)$ is the bounded derived category of a certain modular category. Our goal is to understand the Hecke subcategory $e_\L\D_G(G)$. As a first step we will study the monoidal category $e_\L\D_U(G)$ and some structures on it.

In \S\ref{torus}, we prove the key result that there are no non-trivial idempotents in $\D(T)$, where $T$ is a torus. We will use this result to show that the idempotent associated with an admissible pair on a solvable group is a minimal idempotent.

In \S\ref{solvable}, we study the case of solvable groups. In \S\ref{pfmain1} we will prove
\bthm\label{main1}
Let $G=UT$ be a connected solvable group with unipotent radical $U$ and a maximal torus $T$. Let $(H,\L)$ be a Heisenberg admissible pair for $G$.  Then there exists a canonical equivalence of monoidal categories $e_\L\D_U(G)\cong e_\L\D_U(U)\boxtimes \D(T)$.
\ethm
Using this in \S\ref{pfmain2} we will prove
\bthm\label{main2}
Let $G$ be a group such that $G^\circ$ is solvable. Let $(H,\L)$ be an admissible pair for $G$ with normalizer $G'$. Then $e_\L\in \DGp$ is a minimal idempotent and $f_{\L}\in \DG$ is also a minimal idempotent.
\ethm

In \S\ref{pfmain3} we prove that
\bthm\label{main3}
Let $G$ be such that $G^\circ$ is solvable. Then Conjecture \ref{conjmain} holds for $G$.
\ethm 
We prove this result by proving the following:
\bthm\label{main4}
Let $G$ be a group such that for any closed subgroup $G_1\subset G$ and any admissible pair $(H,\L)$ for $G_1$, $f_\L\in \D_{G_1}(G_1)$ is a minimal idempotent. 
\bit
\item[(i)] Suppose $f\in \DG$ is a minimal idempotent. Then there exists an admissible pair $(H,\L)$ for $G$ such that $f_\L\cong f.$
\item[(ii)] Suppose $(H,\L)$ is an admissible pair for $G$ with normalizer $G'.$ Then $f_\L$ is a locally closed idempotent in $\DG$ and we have a braided equivalence 
\beq
\indg : e_\L\DGp\to f_\L\DG.
\eeq
\eit
\ethm

\brk
As a corollary of the previous result, we see that Conjecture \ref{conjweak} implies Conjecture \ref{conjmain}.
\erk

Finally in \S\ref{redtoheis}, we briefly sketch a possible approach towards eventually developing a theory of character sheaves on general affine algebraic groups. In particular, we observe that our results reduce the problem of defining and studying character sheaves on general groups to the Heisenberg case.

\section{Averaging and induction functors}\label{avgindfun}
In this section, we recall  the definition and properties of the averaging and induction functors from \cite[\S5,\S6]{BD08}\footnote{The results from {\it op. cit.} are proved for unipotent groups. We modify some of the proofs from  {\it op. cit.} and \cite{B} to prove these results in general.}. First we define the averaging functor and use it to define the induction functor.

\subsection{Averaging functor}
Let $G$ be a perfect quasi-algebraic group acting on a perfect quasi-algebraic scheme $X$ and let $G'\subset G$ be a closed subgroup. Then we have a ``forgetful functor'' $\D_G(X)\to \D_{G'}(X)$. Again, by a slight abuse of notation, we will often use the same symbol to denote an object of $\D_G(X)$ and its forgetful image in $\D_{G'}(X)$. We now describe the averaging with compact supports functor $\av_{G/G'}:\D_{G'}(X)\to \D_G(X).$ Informally, for $M\in \D_{G'}(X)$, $\av_{G/G'}(M)$ is obtained by looking at all conjugates $g^*M$ of $M$ for $G\ni g:X\to X$ and ``summing up over the quotient $G/G'$''. More precisely, we have an equivalence $\D_G(G/G'\times X)\cong \D_{G'}(X)$ where $G$ acts on $G/G'\times X$ diagonally. Let us denote this equivalence by $\D_{G'}(X)\ni M\mapsto \t{M}\in \D_G(G/G'\times X).$ Then we define averaging with compact support as the composition $$\D_{G'}(X)\cong\D_G(G/G'\times X)\stackrel{{p_X}_!}{\longrightarrow}\D_G(X)$$ i.e. $\av_{G/G'}(M):={p_X}_!(\t{M})$, where $p_X:G/G' \times X\to X$ is the second projection. 

In the language of stacks, we have a morphism of stacks $X/G'\to X/G$ and the forgetful functor and averaging functor are the pullback and pushforward functors associated with this morphism.

Suppose we have an intermediate closed subgroup $G'\subset K\subset G$. Consider the diagram
$$\xymatrix{
K/G'\times X \ar@{^{(}->}[r]\ar[d] & G/G'\times X\ar[d]\\
X\ar@{^{(}->}[r]^-{i} & G/K\times X\ar[d]\\
 & X\\
}$$
where the upper square is a pullback square. Note that the equivalence $\D_G(G/K\times X)\stackrel{\cong}{\longrightarrow}\D_K(X)$ is given by  (the forgetful functor followed by) the pullback $i^*$.
Using this,  we obtain a natural isomorphism of functors
\beq\label{comofavg}
\av_{G/K}\circ \av_{K/G'}\stackrel{\cong}{\longrightarrow} \avg.
\eeq
Also we see that for each $M\in \D_{G'}(X)$, we have a natural morphism $\avg M\to \av_{K/G'}M$ in $\D_K(X)$ where we consider $\avg M$ as an object of $\D_K(X)$ via the forgetful functor. Informally, the functor $\av_{K/G'}$ is obtained by ``summing up'' over $K/G'$ which is a closed subscheme of $G/G'$ and the natural morphism is obtained by a push forward with compact support (${p_X}_!$) of the adjunction $\id_{\D_K(G/G'\times X)}\to i_!i^*$ where $i$ denotes the inclusion $K/G'\times X\hookrightarrow G/G'\times X$. In particular taking $G'=K$ we obtain
\blem\label{forgetful}
For each $M\in \D_{G'}(X)$, there is a functorial morphism $$\avg M\to M$$ in $\D_{G'}(X)$.
\elem

By \cite[Defn. 6.3]{BD08}, we have functorial isomorphisms
\beq
L\otimes \avg(M)\stackrel{\cong}{\longrightarrow}\avg(L\otimes M), \hbox{ for } L\in\D_G(X), M\in\D_{G'}(X),
\eeq
where on the right hand side, $L$ is considered as an object of $\D_{G'}(X)$ via the forgetful functor.

Suppose now that $G$ acts on $X,Y$ and that $f:X\to Y$ is a $G$-morphism. Then it is easy to see that we have a natural isomorphisms of functors (see \cite[Defns. 6.1,6.2]{BD08})
\beq\label{natisomavg}
f_!\circ\avg\stackrel{\cong}{\longrightarrow} \avg\circ f_!: \D_{G'}(X)\to \D_G(Y),
\eeq
\beq
f^*\circ\avg\stackrel{\cong}{\longrightarrow} \avg\circ f^*: \D_{G'}(Y)\to \D_G(X).
\eeq

Now let $X,Y$ be as above and let $G$ act on $X\times Y$ diagonally. For $M\in\D_{G'}(X), N\in\D_{G'}(Y)$  we will construct functorial morphisms 
\beq\label{funmor}
(\av_{G/G'}M)\boxtimes (\avg N)\to \avg(M\boxtimes N).
\eeq

We use the following lemma, which is an immediate consequence of proper base change and the projection formula:
\blem\label{prodofsums}
Suppose we have a pullback square of perfect quasi-algebraic schemes
$$\xymatrix{
X \ar[r]^{\t{f_1}}\ar[d]_{\t{f_2}}\ar[rd]^f & X_2\ar[d]^{f_2}\\
X_1\ar[r]_{f_1} & Y.\\
}$$
Then for $M_1\in\D(X_1)$ and $M_2\in \D(X_2)$, we have functorial isomorphisms $$f_!(\t{f_2}^*M_1\otimes \t{f_1}^*M_2)\cong {f_1}_!M_1\otimes {f_2}_!M_2.$$
\elem

Now coming back to the situation above, with $G$ acting on $X,Y$ we have the following diagram of 3 pullback squares
$$\xymatrix{
G/G'\times X\ar[d]_{p_X} & G/G'\times X\times Y\ar[l]\ar[d] & G/G'\times X\times G/G'\times Y\ar[l]\ar[d]\ar[ld]^{P_{X,Y}}\\
X & X\times Y\ar[l]\ar[d] &  X\times G/G'\times Y\ar[l]\ar[d]\\
& Y & G/G'\times Y,\ar[l]^{p_Y}\\
}
$$
where all the maps are projections onto the obvious coordinates. Now using proper base change and Lemma \ref{prodofsums}, we have natural isomorphisms 
\beq\label{tenofav} 
(\av_{G/G'}M)\boxtimes (\avg N)={p_X}_!\t{M}\boxtimes {p_Y}_!\t{N}\cong {P_{X,Y}}_!(\t{M}\boxtimes \t{N}).
\eeq
 On the other hand \beq \avg(M\boxtimes N)={p_{X\times Y}}_!(\t{M\boxtimes N}), \eeq where $\t{M\boxtimes N}\in \D_G(G/G'\times X\times Y).$

Consider the diagram
$$\xymatrixcolsep{4pc}\xymatrix{
G/G'\times X\times Y\ar@{^{(}->}[d]_{\Delta}\ar[rd]^{p_{X\times Y}} & \\
G/G'\times X\times G/G'\times Y\ar[r]_-{P_{X,Y}} & X\times Y,\\
}
$$
where $\Delta(gG',x,y):=(gG',x,gG',y)$. We can check that we have a natural isomorphism $$\t{M\boxtimes N}\cong \Delta^*(\t{M}\boxtimes \t{N})$$ and hence 
\beq\label{avoften} 
\avg(M\boxtimes N)={P_{X,Y}}_!\Delta_{!}\Delta^*(\t{M}\boxtimes \t{N}). 
\eeq

Now we have the natural adjunction morphism $\t{M}\boxtimes \t{N}\to\Delta_*\Delta^*(\t{M}\boxtimes \t{N})$. Since $\Delta_*\cong\Delta_!$, by applying ${P_{X,Y}}_!$ to the adjunction morphism above, we obtain the functorial morphism (\ref{funmor}).

Similarly in the situation above, we can prove the following:
\bprop\label{avgasbimod} (See \cite[Prop. 6.7]{BD08}.) 
For  $M\in \D_{G'}(X)$ and $N\in \D_G(Y)$ there exist functorial isomorphisms 
\beq
(\avg M)\boxtimes N \stackrel{\cong}{\longrightarrow} \avg (M\boxtimes N).
\eeq
\eprop

\subsection{Averaging for actions on groups by automorphisms}

Now suppose that $G$ acts on a perfect quasi-algebraic group $H$ by group automorphisms. In this case, $\D_G(H)$ and $\D_{G'}(H)$ are monoidal categories under convolution with compact support and the forgetful functor is monoidal. Now we prove that
\bprop
For $M,N\in \D_{G'}(H)$, there exist functorial morphisms
\beq
\phi(M,N):(\avg M)\ast (\avg N)\to \avg(M\ast N)
\eeq which provide the functor $\avg:\D_{G'}(H)\to \D_{G}(H)$ with a \emph{weak semigroupal structure} (see \cite[Def. 2.4]{BD08}).
\eprop
\bpf
By (\ref{funmor}), we have a functorial morphism $(\av_{G/G'}M)\boxtimes (\avg N)\to \avg(M\boxtimes N)$. Let $\mu_H:H\times H\to H$ be the multiplication on $H$. Then we take the pushforward (with compact supports) of the previous morphism along $\mu_H$ and $\phi(M,N)$ is defined as the composition 
\beq{\label{phimn}}
{\mu_H}_!\left((\av_{G/G'}M)\boxtimes (\avg N)\right)\to {\mu_H}_!\avg(M\boxtimes N)\cong \avg{\mu_H}_!(M\boxtimes N), 
\eeq
where the last isomorphism is the one provided in (\ref{natisomavg}). It is straighforward to verify that this provides a weak semigroupal structure to the functor $\avg:\D_{G'}(H)\to \D_{G}(H).$
\epf

The morphisms $\phi(M,N)$ are not isomorphisms in general. The next proposition, which gives a criterion under which $\phi(M,N)$ is an isomorphism, will be used crucially in this text. For $N\in \D_{G'}(H)$ and $g\in G$ we let ${ }^{g}N:={(g^{-1})^*N}$ where $g^{-1}:H\to H$ is the action of $g^{-1}\in G$ on $H$.
\bprop\label{precri}
Suppose $M,N\in \D_{G'}(H)$ are such that $M\ast { }^{g}N=0$ for each $g\in G-G'$. Then the morphism 
$\phi(M,N):(\avg M)\ast (\avg N)\stackrel{\cong}{\longrightarrow} \avg(M\ast N)$ is an isomorphism.
\eprop
\bpf
By (\ref{phimn}), it is enough to check that $${\mu_H}_!\left((\av_{G/G'}M)\boxtimes (\avg N)\right)\to {\mu_H}_!\avg(M\boxtimes N)$$ is an isomorphism. Setting $X=Y=H$ in (\ref{tenofav}) and (\ref{avoften}), we must show that the morphism ${\mu_H}_!{P_{H,H}}_!(\t{M}\boxtimes \t{N})\to{\mu_H}_!{P_{H,H}}_!\Delta_{!}\Delta^*(\t{M}\boxtimes \t{N})$  obtained using the adjunction $\t{M}\boxtimes \t{N}\to \Delta_{!}\Delta^*(\t{M}\boxtimes \t{N})$ is an isomorphism.

Note that for $g\in G$, the restriction $\t{M}|_{\{gG'\}\times H}$ is isomorphic to ${ }^{g}M$ after identifying $\{gG'\}\times H$ with $H$. We have the following commutative diagram
$$\xymatrixcolsep{4pc}\xymatrix{
\{g_1G'\}\times H \times\{g_2G'\}\times H\ar@{^{(}->}[r]\ar[d]_{\mu_H} & G/G'\times H\times G/G'\times H\ar[r]^-{P_{H,H}}\ar[d]_{\id\times \mu_H} & H\times H\ar[d]^{\mu_H}\\
\{g_1G'\} \times\{g_2G'\}\times H\ar@{^{(}->}[r] & G/G'\times G/G'\times H\ar[r] & H.\\
}
$$
The restriction of $(\id\times \mu_H)_!(\t{M}\boxtimes \t{N})$ to $\{g_1G'\} \times\{g_2G'\}\times H$ can be identified with ${ }^{g_1}M\ast { }^{g_2}N={ }^{g_1}(M\ast { }^{g_1^{-1}g_2}N).$ By the given criterion, this is zero unless $g_1G'=g_2G'.$ Hence we see that the morphism $(\id\times \mu_H)_!(\t{M}\boxtimes \t{N})\to (\id\times \mu_H)_!\Delta_!\Delta^*(\t{M}\boxtimes \t{N})$ is an isomorphism. Hence from the commutative diagram above, we get our desired isomorphism.
\epf

Using Proposition \ref{avgasbimod} we can prove:
\bprop\label{avgconvbimod} (See \cite[Cor. 6.8]{BD08}.)
Let $H$ be equipped with an action of $G$ as above. Then for $M\in \D_{G'}(H)$ and $N\in \D_G(H)$ there exist functorial isomorphisms 
\beq
(\avg M)\ast N \stackrel{\cong}{\longrightarrow} \avg(M\ast N),
\eeq
\beq
N\ast (\avg M)\stackrel{\cong}{\longrightarrow} \avg(N\ast M).
\eeq

\eprop

Recall (from Lemma \ref{forgetful}) that we have a functorial morphism $\avg M\to M$ for each $M\in \D_{G'}(H)$. Hence we have a functorial morphism  $L\ast(\avg M)\to L\ast M$ in $\D(H)$ for each $L\in \D(H)$. Using similar methods as above, we can prove
\bprop\label{LastM} (See also \cite[Prop. 6.11]{BD08}.)
Let $L\in\D(H)$ and $M\in\D_{G'}(H)$ be such that $L\ast { }^{g}M=0$ for all $g\in G-G'.$ Then the functorial morphism $$L\ast(\avg M)\stackrel{\cong}{\longrightarrow} L\ast M$$ in $\D(H)$ is an isomorphism.
\eprop

\subsection{Fourier-Deligne transform and averaging}
We recall some facts about the equivariant Fourier-Deligne transform and averaging from \cite[\S3.4,\S6.1]{BD08}. Let $G$ act on $H$ as before, where we now assume that $H$ is (perfect) unipotent and connected. Let $H^*$ be the Serre dual of $H$ and let $\sE$ be the universal local system (which has a canonical $G$-equivariant structure) on $H\times H^*$. Let $\pr:H\times H^*\to H$ and $\pr':H\times H^*\to H^*$ be the two projections.
\bdefn
We define the equivariant Fourier-Deligne transform functors $\FDG,\iFDG$ by
\beq
\FDG:\D_G(H)\to \D_G(H^*), \FDG(N):=\can_H\otimes \pr'_!(\sE^{-1}\otimes\pr^*N),
\eeq
\beq
\iFDG:\D_G(H^*)\to \D_G(H), \iFDG(M):=\can_H\otimes \pr_!(\sE\otimes\pr'^*M),
\eeq
\edefn

\bprop\label{FDsemigroupal} (See \cite[\S3.5]{BD08}.)
We have natural isomorphisms
\beq
\iFDG(M_1)\ast \iFDG(M_2)\stackrel{\cong}{\longrightarrow} \iFDG(M_1\otimes M_2) \hbox{ for }M_1,M_2\in \D_G(H^*)
\eeq
which provides the functor 
\beq
\iFDG:(\D_G(H^*),\otimes)\to (\D_G(H),\ast)
\eeq
with a strong semigroupal structure. In case $H$ is \emph{commutative}, $\iFDG$ is a monoidal equivalence.
\eprop

Let $K\subset [H,H]^*$ denote the image of the restriction homomorphism $H^*\to [H,H]^*$. For $\E\in K$, we have the closed idempotent $e_\E:=\E\otimes \can_{[H,H]}$ in $\D(H)$ (extended by zero outside $[H,H]$). Note that any $\E\in K$ is an $H$-equivariant multiplicative local system on $[H,H]$ for the conjugation action. Since $H$ is connected unipotent, we can consider $\D_H(H)$ as a full subcategory of $\D(H)$. We can check that $e_\E\D(H)=\D_{[H,H],\E}(H)\subset \D_H(H)$ where we consider all categories as full subcategories of $\D(H)$.
\bthm\label{FDmain}(See \cite[Thm. 3.16]{BD08}.) Using the above notation, we have
\bit
\item[(i)] For $\E\in K$, the restriction of $\iFD$ induces a monoidal equivalence of monoidal categories 
\beq
\iFD|_{\D(H^*_\E)}:(\D(H^*_\E),\otimes)\stackrel{\cong}{\longrightarrow} e_\E\D(H)\cong e_\E\D_H(H).
\eeq
We have $\iFD((\Qlcl)_{H^*_\E})\cong e_\E$, a closed idempotent in $\D(H)$, where $(\Qlcl)_{H^*_\E}$ is the constant sheaf supported on the connected component $H^*_\E$. 
\item[(ii)] The Fourier-Deligne transform of the constant sheaf on $H^*$ is given by \beq \iFD((\Qlcl)_{H^*})\cong \bigoplus\limits_{\E\in K}e_\E=:e_0.\eeq
\item[(iii)] The object $e_0\in \D(H)$ is a closed idempotent.
\item[(iv)] The Fourier-Deligne transform defines a monoidal equivalence
\beq
\iFD: (\D(H^*),\otimes)\stackrel{\cong}{\longrightarrow} e_0\D(H)\cong e_0\D_H(H).
\eeq
\item[(v)] The equivariant Fourier-Deligne transform defines a monoidal equivalence
\beq
\iFDG: (\D_G(H^*),\otimes)\stackrel{\cong}{\longrightarrow} e_0\D_G(H).
\eeq
\eit
\ethm

The equivariant Fourier-Deligne transform is compatible with averaging:
\bprop\label{FDandavg}(See \cite[Prop. 6.4]{BD08}.)
For $G'\subset G$ and $H$ as before,  we have a natural isomorphism of functors
\beq
\iFDG\circ\avg \stackrel{\cong}{\longrightarrow} \avg\circ\iFD^{G'}:\D_{G'}(H^*)\to \D_G(H).
\eeq
\eprop

\subsection{Induction functor}\label{indfunc}
As before, let $G'$ be a subgroup of $G$. The inclusion $G'\stackrel{i}{\hookrightarrow} G$ is a $G'$-morphism for the conjugation action of $G'$ on $G$ and $G'$. We have the restriction functor $\Resg:\DG\to \DGp$ defined by composing $i^*$ with the forgetful functor $\D_G(G)\to \D_{G'}(G)$.
\bdefn
The induction with compact supports is defined as the composition of functors 
\beq
\ind_{G'}^G:\D_{G'}(G')\stackrel{i_!}{\longrightarrow}\D_{G'}(G)\stackrel{\avg}{\longrightarrow}\D_G(G).
\eeq
The functor $\ind_{G'}^G$ has the structure of a weak semigroupal functor coming from the monoidal structure of $i_!$ and the weak semigroupal structure of $\avg$.
\edefn

From (\ref{comofavg}), we see that
\blem\label{compind}
If we have closed subgroups $G'\subset K\subset G$, then we have a natural isomorphism of functors
\beq
\ind_K^G\circ\ind_{G'}^K\cong \indg.
\eeq
\elem

Using Lemma \ref{forgetful}, we get
\blem
For each $M\in \DGp$ we have a functorial morphism $$\indg M\to i_!M$$ in $\D_{G'}(G)$ where we consider both objects in the morphism above as objects in $\D_{G'}(G)$, and from the above we obtain a functorial morphism $$\Resg\circ\indg M\to M$$ in $\D_{G'}(G')$.
\elem

\bprop\label{macrit}
Suppose $M,N\in\D_{G'}(G')$ are such that $M\ast\delta_g\ast N=0$ for each $g\in G-G'$, where $\delta_g$ is the delta-sheaf on $G$ supported at $g$ and where we consider $M,N$ as sheaves on $G$ by extending by zero outside $G'.$ Then the weak semigroupal structure morphism is an isomorphism:
\beq
\phi(M,N): (\indg M)\ast(\indg N)\stackrel{\cong}{\longrightarrow}\indg(M\ast N).
\eeq
\eprop
\bpf
This follows from Proposition \ref{precri}, since the condition $M\ast\delta_g\ast N=0$ is equivalent to the condition $M\ast\delta_g\ast N\ast \delta_{g^{-1}} = M\ast { }^gN=0$.
\epf

The following result is proved in \cite{BD08}.
\bprop\label{combra} (See \cite[Cor. B.47]{BD08}.)
The weak semigroupal structure on $\indg$ is compatible with the braidings in the categories $\D_{G'}(G')$ and $\D_G(G)$.
\eprop

\section{Idempotents in $\D(G)$ and induction}\label{idandind}
Let $G'\subset G$ be a closed subgroup. In this section, we study idempotents in $\D_{G'}(G')$ which satisfy the  Mackey condition (see below) with respect to $G$ and study induction of such idempotents. Most of these results are proved in \cite[\S5]{BD08} and \cite{B} only for (perfect) unipotent groups, nevertheless many of the proofs also work for all (perfect quasi-) algebraic groups, perhaps after minor adjustments. 
\subsection{Idempotents satisfying the Mackey condition}
\bdefn
An idempotent $e\in\D_{G'}(G')$  is said to satisfy the {\it Mackey condition} with respect to $G$ if for every $x\in G-G'$ we have $e\ast \delta_x\ast e = 0$, where we view $e$ as an object of $\D(G)$ by extension by zero outside $G'$ and where $\delta_x$ is the delta-sheaf supported at the point $x$.
\edefn

\brk
Let $(H,\L)\in\P(G)$ be a pair with normalizer $G'$. We have the corresponding closed idempotent $e_\L:=\L\otimes \can_H \in \D_{G'}(G').$ Then the above Mackey condition for the idempotent $e_\L$ is equivalent to the Mackey condition(iii) in Definition \ref{defadm}
\erk

\subsection{Induction of idempotents satisfying the Mackey condition}
Using  Propositions \ref{avgconvbimod}, \ref{LastM} and \ref{macrit} we can prove
\bprop\label{weakid} (See \cite[\S5]{BD08},\cite[\S5.8]{B}.)
Let $e\in \D_{G'}(G')$ be an idempotent satisfying the Mackey condition with respect to $G$. Then 
\bit
\item[(i)] The object $f=\indg e$ is an idempotent in $\D_G(G)$. 
\item[(ii)] If $M\in e\D_{G'}(G')$, then $\indg M\in f\DG.$
\item[(iii)] If $N\in \DG$, then $(e\ast N)|_{G-G'}\cong (e\ast N\ast e)|_{G-G'}=0.$ Hence $e\ast N\cong e\ast (N|_{G'})$ and by a slight abuse of notation, we can consider $e\ast N$ as an object of $e\DGp$.
\item[(iv)] For each $N\in \DG$, there is an isomorphism
\beq
f\ast N \stackrel{\cong}{\longrightarrow} \indg (e\ast N)
\eeq
functorial with respect to $N$.
\item[(v)] For each $M\in e\DGp$ we obtain functorial isomorphisms 
\beq
e\ast (\indg M) \stackrel{\cong}{\longrightarrow} e\ast M \hbox{ in } \D_{G'}(G),
\eeq
\beq
e\ast (\Resg\circ\indg M) \stackrel{\cong}{\longrightarrow} e\ast M \hbox{ in } \DGp,
\eeq
by taking the convolution of $e$ with the canonical morphism $\indg M\to M$ in $\D_{G'}(G)$ and $\Resg\circ\indg M\to M$ in $\DGp$. In particular we see that $e\ast f\cong e.$

\item[(vi)] The restriction \beq\label{indfun} \indg|_{e\DGp}:e\DGp\to f\DG \eeq is strong braided semigroupal and induces a bijection on isomorphism classes of objects.

\item[(vii)] The strong semigroupal functor (\ref{indfun}) induces a bijection between the sets of isomorphism classes of idempotents  in $e\DGp$ and $f\DG$. In particular, if $e$ is minimal in $\DGp$ then $f$ is minimal in $\DG$.
\item[(viii)] Let $e\neq 0$. Then $f\neq 0$. Let $(A,\E)\in \P(G)$ be such that $f\ast e_\E\cong f$, where $e_\E=\E\otimes \can_A$ is an idempotent in $\D(G)$. Then we must have $e\ast e_\E\cong e$ and $A\subset G'$.
\item[(ix)] If the functor $M\mapsto e\ast M$ is isomorphic to the identity functor on $e\DGp$, the functor (\ref{indfun}) is faithful.
\item[(x)] If the functors $M\mapsto e\ast M$ and $N\mapsto f\ast N$ are isomorphic to the identity functors on $e\DGp$ and $f\DG$ respectively, the functor (\ref{indfun}) is a braided monoidal equivalence of braided monoidal categories, a quasi-inverse to which is provided by the functor (also see (iii) above) \beq f\DG\ni N\mapsto e\ast N \in e\DGp.\eeq
\eit
\eprop

\subsection{Locally closed idempotents via averaging}
Suppose $G$ is a perfect quasi-algebraic group acting by group automorphisms on a perfect connected unipotent group $H$. We will construct certain locally closed idempotents in $\D_G(H)$ using averaging functors. 
\brk\label{lcremark}
This construction was carried out in the case of unipotent  $G$ in \cite[Prop. 7.17]{BD08}. In this special case, the idempotents that are obtained by averaging in fact turn out to be closed idempotents. The reason is that the orbits of a unipotent group acting on an affine scheme are always closed. For general groups, we can only say that the orbits are locally closed.
\erk
\bprop\label{locallyclosed}
Let $\L$ be a multiplicative local system on $H$ and let $G'\subset G$ denote the stabilizer of $\L$ for the action of $G$ on $H^*.$ Then $e_\L:=\L\otimes \can_H$ is a closed idempotent in $\D_{G'}(H)$ and $\avg(e_\L)$ is a locally closed idempotent in $\D_G(H)$.
\eprop
\bpf
If $g\in G-G'$, then $\L\ncong { }^g\L$ and we see that $e_\L\ast { }^ge_\L=0$. Hence by Proposition \ref{precri} $\avg(e_\L)$ is an idempotent in $\D_G(H).$ We will now use the equivariant Fourier-Deligne transform. Let $\delta_\L\in \D_{G'}(H^*)$ denote the delta-sheaf supported at $\L\in H^*.$ Then we have $\iFD^{G'}(\delta_\L)\cong e_\L.$ Hence by Proposition \ref{FDandavg} \beq \avg(e_\L)\cong \avg(\iFD^{G'}(\delta_\L))\cong \iFDG(\avg(\delta_\L)). \eeq
Let $G\L\subset H^*$ denote the $G$-orbit of $\L$ in $H^*$. Then $G\L\subset H^*$ is locally closed and $\avg(\delta_\L)\cong (\Qlcl)_{G\L}\in \D_G(H^*)$, where $(\Qlcl)_{G\L}$ is the constant sheaf on the orbit $G\L$ extended by zero outside $G\L$. Since the orbit is locally closed, we see that $(\Qlcl)_{G\L}\in \D_G(H^*)$ is a locally closed idempotent. Hence by Theorem \ref{FDmain}(v), $\avg(e_\L)\cong \iFDG\left((\Qlcl)_{G\L}\right)$ is a locally closed idempotent in $e_0\D_G(H)$ where $e_0$ is as before. By Lemma \ref{lcinlc} we conclude that  $\avg(e_\L)$ is a locally closed idempotent in $\D_G(H)$.
\epf

Using Propositions \ref{weakid} and \ref{locallyclosed} we can prove
\bcor\label{loccloscor} (See \cite[Lem. 7.18]{BD08}.)
Let $(A,\E)\in \Pn(G)$ be a pair with normalizer $G_1\subset G$. Consider the corresponding closed idempotent  $e_\E=\E\otimes\can_A\in \D_{G_1}(G_1)$. We have
\bit
\item[(i)] $e_\E\in \D_{G_1}(G_1)$ satisfies the Mackey condition with respect to $G$.
\item[(ii)] Let $f_\E:=\ind_{G_1}^G(e_\E)$. Then $f_\E$ is a locally closed idempotent in $\DG$ and hence $f_\E\DG$ is a braided monoidal category.
\item[(iii)] We have an equivalence of braided monoidal categories
\beq
\ind_{G_1}^G:e_\E\D_{G_1}(G_1)\to f_\E\DG.
\eeq
\item[(iv)] $e_\E\ast f_\E\cong e_\E$.
\eit
\ecor

\section{Analysis of admissible pairs}\label{anadmpair}

\subsection{Admissible pairs on a reductive group}

Let us first prove that there are no non-trivial admissible pairs on reductive groups.
\bprop\label{admonred}
If $(H,\L)$ is an admissible pair on a reductive group $G$, then $H=\{1\}$.
\eprop
\bpf
The connected unipotent subgroup $H$ lies inside a maximal unipotent $U\subset G$. Now there exists a $g\in G$ such that $U\cap { }^gU=\{1\}$. If $H$ is of positive dimension, then this $g$ violates condition (iii) of the definition of admissible pair.
\epf

\subsection{Extensions of Heisenberg admissible pairs}\label{eoap}
Let $(H,\L)$ be a Heisenberg admissible pair on the group $G$. This means that $H\normal G$ and $\L$ is $G$-equivariant for the conjugation action on $H$. Let $U:=\R_u(G)$. We have an action of the reductive group $\Gamma=G/U$ on the commutative unipotent group $U/H$ and the isogeny $\phi_\L: U/H\to (U/H)^*$ is $\Gamma$-equivariant. 

In this subsection, we will classify all admissible pairs $(L,\L')$ extending $(H,\L)$ i.e. we will determine when a pair $(L,\L')\geq (H,\L)$ is also admissible.

Suppose first that  $(H,\L)\leq (L,\L')$ and that $H\subset L\subset U$. Let us restrict the commutator map $c:U\times U\to H$ to $c_L:L\times L\to H$. The skew-symmetric biextension $\phi_\L$ is obtained by pulling back the local system $\L$ on $H$ along the commutator map. By \cite[Prop. 7.7.]{B}, we must have that $c_L^*\L$ is the trivial local system on $L\times L.$ This is equivalent to the condition that $L/H\subset U/H$ is an isotropic (see Appendix \ref{equivisotropic}) subgroup for the skew-symmetric isogeny $\phi_\L$. 

On the other hand, suppose $H\subset L\subset  U$ is such that $L/H\subset U/H$ is an isotropic subgroup. Hence $c_L^*\L$ is trivial and by \cite[Prop. 7.7.]{B}, there exists a multipicative local system $\L'$ on $L$ such that $\L'|_H\cong \L$. Let $L^*_{\L}\subset L^*$ be the perfect subscheme of multiplicative local systems on $L$ whose restriction to $H$ is isomorphic to $\L$. We have the connected commutative subgroup $(L/H)^*\subset L^*$ of multiplicative local systems on $L$ which are trivial when restricted to $H$ and $L^*_\L\subset L^*$ is a coset for this subgroup. 
	
Since $\L$ is $U$-equivariant, the conjugation action of $U$ on $L$ induces an action of $U$ on $L^*_\L.$ Let  $L\subset L^\perp\subset U$ be such that $(L/H)^\perp=L^\perp/H$ with respect to the skew-symmetric biextension $\phi_\L$, or in other words $L^\perp$ is the kernel of the composition $U\onto U/H\stackrel{\phi_\L}{\onto}(U/H)^*\onto (L/H)^*$. This composition is given by $u\mapsto c_u|_L^*\L$, where $c_u|_L:L\to H$ is the commutator with $u$ map. We can check that:
\blem\label{transitivity}
The action of $U$ on $L^*_\L$ is given by ${ }^u\L'=c_u|_L^*\L\otimes \L'$ for $u\in U, \L'\in L^*_\L$. Hence this action is transitive and for any $\L'\in L^*_\L$, $\N_U(L,\L')=\Stab_U(\L')=L^\perp.$
\elem

\bcor\label{compatibility}
Suppose $H\subset L\subset U$ and $(H,\L)\leq (L,\L')$. Then for every $g\in G$, there exists $u\in U$ such that the pairs $(L,\L')$ and $({ }^{ug}L,{ }^{ug}\L')$ are compatible.
\ecor
\bpf
Since the pair $(H,\L)$ is normalized by all of $G$,  $(H,\L)\leq({ }^{g}L,{ }^{g}\L')$. Now by \cite[Prop. 7.7.]{B} we can extend ${ }^{g}\L'|_{(L\cap { }^{g}L)^\circ}$ to a multiplicative local system $\L''$ on $L$. Hence $(L,\L'')$ is compatible with the pair $({ }^{g}L,{ }^{g}\L')$. It is clear that $\L''\in L^*_\L$ and hence by Lemma \ref{transitivity}, there exists $u$ such that ${ }^{u}\L''\cong \L'$. This is the required $u$.
\epf

\bprop\label{Lnormal}
Let $(L,\L')\in \P(G)$ be any pair such that $(H,\L)\leq (L,\L')$. (Here we do not assume $L\subset U$.) Then $(L,\L')$ is an admissible pair for $G$ if and only if $L\normal G$. 
\eprop
\bpf
First suppose that $L\normal G$. Hence in this case $L$ lies inside the unipotent radical $U$. By Corollary \ref{compatibility}, for each $g\in G$, there exists a $u\in U$ such that $ug\in \N_G(L,\L')$. Hence by Lemma \ref{transitivity} we have an extension
\beq\label{normalizer}
0\to L^\perp\to \N_G(L,\L')\to \Gamma\to 0.
\eeq

From (\ref{normalizer}), we see that $\R_u(\N_G(L,\L'))=L^{\perp\circ}$ and clearly $L^{\perp\circ}/L$ is commutative. The skew-symmetric biextension $\phi_{\L'}: L^{\perp\circ}/L \to (L^{\perp\circ}/L)^*$ associated to the pair $(L,\L')$ is the one induced from $\phi_{\L}$ on the subquotient $L^{\perp\circ}/L$ of $U/H$ in the standard way. Hence it follows that $\phi_{\L'}$ is also an isogeny. Since $L\normal G$, the Mackey condition is automatically satisfied. Hence $(L,\L')$ is an admissible pair for $G$.

Now suppose that $L\subset U$ but there is a $g\in G$ which does not normalize $L$. Then by Corollary \ref{compatibility}, there exists $u\in U$ such that  $(L,\L')$ and $({ }^{ug}L,{ }^{ug}\L')$ are compatible. This means that the pair $(L,\L')$ does not satisfy the Mackey condition and hence  cannot be an admissible pair.

Now the only remaining case is when $L$ is not entirely contained inside $U$. Let $U_{max}$ be a maximal connected unipotent subgroup of $G$ containing $L$. Then we know that there exists $g\in G$ such that $U_{max}\cap { }^{g}U_{max}=U$ and hence $U_{max}\cap { }^{ug}U_{max}=U$ for all $u\in U.$ Hence it follows that there is a $g\in G$ such that $L\cap { }^{ug}L\subset U$ for each $u\in U.$ Since we have assumed that $L$ is not contained in $U$, it follows that $ug\notin \N_G(L)\supset \N_G(L,\L')$ for each $u\in U.$ Now, we apply Corollary \ref{compatibility} for this $g$ and  the pair $((L\cap U)^\circ, \L'|_{(L\cap U)^\circ})$ to get a $u\in U$ such that $ug\in G-\N_G(L,\L')$ violates the Mackey condition for the pair $(L,\L')$. Hence such an $(L,\L')$ cannot be an admissible pair.
\epf

Hence given a {\it Heisenberg} admissible pair $(H,\L)$ for $G$, we have characterized all admissible pairs $(L,\L')$ that extend it. They correspond to $\Gamma$-stable isotropic subgroups $L/H\subset U/H$ and any extension $\L'\in L^*_\L$ of the multiplicative local system $\L$.

\subsection{Analysis of associated idempotents}
Suppose $(H,\L)$ is a Heisenberg admissible pair for $G$ and let $(L,\L')\geq (H,\L)$ be an admissible pair as above. Corresponding to the admissible pairs above we have the idempotents $e_{\L'}=\L'\otimes\can_L\in \DGp$ and $\L\otimes \can_H =e_\L=f_\L, f_{\L'}=\indg(e_{\L'})\in \DG$. Let $G'=\N_G(L,\L').$ In this section, we prove
\bprop\label{heisext}
Using the notation above, we have $e_\L\cong f_{\L'}$ and a braided monoidal equivalence 
\beq
\indg: e_{\L'}\DGp\xto{\cong} e_\L\DG.
\eeq
\eprop
\bpf
By Proposition \ref{weakid}(x), the braided monoidal equivalence follows once we prove that $e_\L\cong \indg(e_{\L'})$ since we know that $e_\L\in \DG$ and $e_{\L'}\in \DGp$ are both closed idempotents. 

By Proposition \ref{Lnormal}, $(L,\L')\in \Pn(G).$ Hence we can think of $f_{\L'}$ as the averaging $\avg(e_{\L}')$ for the conjugation action of $G$ on $L$ (and then extending by zero to all of $G$). As in the proof of  Proposition \ref{locallyclosed}, we will use the equivariant Fourier-Deligne transform. We see that $$\iFD^{G'}:\D_{G'}(L^*)\ni \delta_{\L'}\mapsto e_{\L'}\in \D_{G'}(L),$$ where $\delta_{\L'}$ is the delta sheaf supported at $\L'\in L^*.$ By Lemma \ref{transitivity}, the $G$-orbit of $\L'$, $G\L'\subset L^*$ is equal to $L^*_\L$, hence $\avg(\delta_{\L'})\cong (\Qlcl)_{L^*_\L}\in \D_G(L^*)$, the extension by zero of the constant sheaf supported on the orbit $L^*_\L$. By the compatibility of averaging with the Fourier-Deligne transform, the proof would be complete once we prove
\blem
We have $\iFDG((\Qlcl)_{L^*_\L})\cong e_\L\in \D_G(L)$.
\elem
We have the projection $\pr:L\times L^*_\L\to L$. By definition, $\iFDG((\Qlcl)_{L^*_\L})=\can_L\otimes \pr_!(\sE|_{L\times L^*_\L})$. First let us show that $\iFDG((\Qlcl)_{L^*_\L})$ is supported inside $H$. We have $L_\L^*=\L'\otimes L^*_{\Qlcl}$, i.e. $L^*_\L$ is a coset of $L^*_{\Qlcl}$ in $L^*$. Now $L/H$ is a commutative group and we have $L^*_{\Qlcl}=(L/H)^*.$ We have the commutative Fourier-Deligne transform $\D((L/H)^*)\to \D(L/H)$. Under this, the constant sheaf on $(L/H)^*$ maps to the delta-sheaf at the origin of $L/H$. Using this we deduce that  $\iFDG((\Qlcl)_{L^*_\L})$ must be zero outside of $H$. Now restrict the projection map further to the map $\pr_H:H\times L^*_\L\to H$. We have $\sE|_{H\times L^*_\L}\cong \pr_H^*\L.$ Hence we see that $\iFDG((\Qlcl)_{L^*_\L})\cong\can_L\otimes {\pr_H}_!(\sE|_{H\times L^*_\L})\cong \can_L\otimes {\pr_H}_!( \pr_H^*\L)\cong \L\otimes \can_L\otimes {\pr_H}_!(\Qlcl)$. Since we know that $\iFDG((\Qlcl)_{L^*_\L})$ is an idempotent, we must have $\iFDG((\Qlcl)_{L^*_\L})\cong e_\L$.
\epf

\subsection{Extensions of admissible pairs}
Let $(H,\L)$ be any admissible pair for $G$ with normalizer $G'$. Now we can immediately classify all admissible pairs $(L,\L')$ extending the given admissible pair.  If  $(H,\L)\leq (L,\L')$ and both are admissible pairs, then by the Mackey condition, we must have $\N_G(L,\L')\subset \N_G(H,\L)=G'$. 

As before, let $U$ be the unipotent radical of $G'$, let $\Gamma=G'/U$ and let $\phi_\L:U/H\to (U/H)^*$ be the skew-symmetric isogeny corresponding to the admissible pair $(H,\L)$.
\bprop\label{extadmpair}
Let $L/H\subset U/H$ be a $\Gamma$-stable isotropic subgroup. For any $\L'\in L^*_\L$, $(L,\L')$ is an admissible pair for $G$. This characterizes all admissible pairs $(L,\L')\geq (H,\L)$. Moreover, the associated idempotents $f_\L$ and $f_{\L'}$ are isomorphic.
\eprop
\bpf
That such a pair $(L,\L')$ is admissible for $G'$ follows from Proposition \ref{Lnormal}, and the Mackey condition for $g$ lying outside $G'$ follows from the Mackey property of $(H,\L)$. The last assertion follows from Proposition \ref{heisext} and Lemma \ref{compind}.
\epf

\subsection{Geometric reduction process for algebraic groups}
In this subsection, we will see that given any non-zero $M\in \D(G)$, there exists an admissible pair $(H,\L)$ for $G$ which is compatible with $M$. This result is proved for (perfect) unipotent groups $G$ in \cite{BD08}, and only very minor adjustments need to be made in the general case. Hence this subsection is essentially a review of \S4 of \cite{BD08}. 

\bprop\label{compext}(See \cite[Cor. 4.9]{BD08}.)
Let $M\in \D(G)$, and let $(A,\E)\in \P(G)$ be compatible with $M$. Let $A_1\supset A$ be a \emph{connected unipotent} subgroup of $G$ such that $[A_1,A_1] \subset A$. If the pullback of $\E$ by the commutator morphism
$A_1\times A_1\to A$ is trivial then there exists a multiplicative local system $\E_1$ on $A_1$ such that $\E_1|_A\cong \E$ and $(A_1,\E_1)$ is compatible with $M$.
\eprop

\bprop\label{maxnormpair}(See \cite[Prop. 4.4]{BD08}.)
Let $M\in \D(G)$ be non-zero. Suppose $(H,\L)\in \Pn(G)$ is maximal among all pairs\footnote{The pair $(\{1\},\Qlcl)\in \Pn(G)$ is compatible with $M$ and hence there must exist a maximal pair $(H,\L)$.} in $\Pn(G)$ that are compatible with $M$.  If the multiplicative local system $\L$ is $G$-equivariant, then $(H,\L)$ is a Heisenberg admissible pair for $G$.
\eprop
We need to show that $\R_u(G)/H$ is commutative and that $\phi_\L$ is an isogeny. As in \cite{BD08}, this can be proved using 
\blem(See \cite[Lem. 4.10 and Lem. 4.11]{BD08}.)
Let $Z\subset \R_u(G)$ be the preimage of the neutral connected component of the center of $\R_u(G)/H$. Let $\phi_\L:\R_u(G)/H\to (Z/H)^*$ be the group morphism obtained by pulling back $\L$ along the commutator map $\R_u(G)\times Z\to H$. Then the morphism $\phi_\L|_{Z/H}:Z/H\to (Z/H)^*$ is an isogeny. In fact, we must have $Z=\R_u(G)$.
\elem
As in \cite{BD08}, this can be proved using Proposition \ref{compext}.

Now using Proposition \ref{maxnormpair} and an inductive argument on the size of $G$, the following is proved in {\it loc. cit.}
\bthm\label{compadmipair} (See \cite[Thm. 4.5]{BD08}.)
Let $M\in \D(G)$ and let $(A,\E)\in \Pn(G)$ be a pair compatible with $M$. Then there exists an admissible pair $(H,\L)\geq (A,\E)$ which is compatible with $M$.
\ethm
Without loss of generality, we may assume that $(A,\E)\in\Pn(G)$ is a maximal pair compatible with $M$. If $\E$ is $G$-equivariant, then by Proposition \ref{maxnormpair}, $(A,\E)$ is already an admissible pair. If not, we can proceed by induction to the normalizer of the pair which is a strictly smaller group.

\section{Analysis of Heisenberg idempotents}\label{anheisid}
Let $(H,\L)$ be a Heisenberg admissible pair for $G$ and $e=e_\L$ the corresponding closed idempotent in $\D_G(G)$. Recall that in this situation, the pair $(H,\L)$ is normalized by the group $G$. Let $U:=\R_u(G)$ be the unipotent radical of $G$. Then $U/H$ is commutative and the map $\phi_\L: (U/H)\to (U/H)^*$ is a (skewsymmetric) isogeny. Let $(K_\L,\theta)$ be the metric group corresponding to this skewsymmetric isogeny (see \cite[Appendix A.10]{B}), where $K_\L=\ker{\phi_\L}$. Note in particular that $e$ can be thought of as a Heisenberg idempotent on the connected unipotent group $U$. Let $H\subset K\subset U$ be such that $K/H=K_\L$. Let $\Gamma$ denote the reductive group $G/U$. Then $\Gamma$ acts on $U/H$ and the isogeny $\phi_\L$ is $\Gamma$-equivariant.

In this section we extend results of \cite{De} where the special case of a finite group $\Gamma$ was considered and study the Hecke subcategory $e\D_U(G)$.


\subsection{The monoidal category $e\D_U(G)$}\label{crossedbraided}
Our goal is to study the Hecke subcategory $e\D_G(G)$,
but let us first consider the category $e\D_U(G)$. Since $U$ is a connected unipotent group, we may consider this as a full subcategory of $\D(G)$. We will describe the subset of $G$ where objects of $e\D_U(G)$ may be supported. For $g\in G$, we will describle the full subcategory $e\D_U(Ug)\subset e\D_U(G)$ of objects supported on the coset $Ug\subset G.$

Since $U$ is normal in $G$, we have the commutator map $c:U\times G\to U$. For $g\in G$, let $c_g:U\to U$ denote the commutator with $g$ map defined by $c_g(u)=ugu^{-1}g^{-1}=u\cdot { }^g(u^{-1})$. For $g,g_1,g_2\in G$ and $u,u_1,u_2\in U$, we have the following relations:
\beq\label{g1g2}
c_{g_1g_2}(u)=c_{g_1}(u)\cdot { }^{g_1}c_{g_2}(u)
\eeq

\beq
c_g(u_1u_2)= { }^{u_1}c_g(u_2)\cdot c_{g}(u_1)
\eeq

Now suppose that $s\in G$ commutes with $g\in G$. Then we see that
\beq\label{centralizer}
c_g({ }^su)=sus^{-1}gsu^{-1}s^{-1}g^{-1}= sugu^{-1}g^{-1}s^{-1}={ }^sc_g(u).
\eeq
In other words, the map $c_g:U\to U$ is equivariant with respect to  $C_G(g)$, the centralizer of $g$ in $G$.

Let $U_g:=\{u\in U|c_g(u)\in H\}$. Then $U_g$ is a subgroup of $U$ containing $H$. Since $H$ is a normal subgroup of $G$, (\ref{centralizer}) implies that $U_g$ (and hence its neutral connected component $U^\circ_g$) is stable under the action of $C_G(g)$ by conjugation. Let $c_g':U_g^\circ\to H$ denote the restriction of the commutator with $g$ map to the neutral connected component of the subgroup $U_g$. Equation (\ref{centralizer}) implies that $c_g'$ is a $C_G(g)$-equivariant map. The subgroup $U_g$ depends only on the coset $Ug$ and ${c'_g}^*\L$ is a multiplicative local system on $U^\circ_g$ which is trivial on $H\subset U^\circ_g.$ Since $\L$ is $G$-equivariant, it follows that $c_g'^*\L$ is $C_G(g)$-equivariant.

We have an action of $G$ on the connected commutative unipotent group $V:=U/H$. We have $U_g/H=V^g\subset V$, the subgroup of $g$-fixed points of $V$. It is clear that $V^g$ and hence $(V^g)^\circ$ are $C_G(g)$-stable. Since $c'^*_g\L$ is trivial on $H$, we may regard $c'^*_g\L$ as a $C_G(g)$-equivariant multiplicative local system on $(V^g)^\circ$.

By (\ref{g1g2}), we see that for $g_1,g_2\in G$, $U_{g_1}\cap U_{g_2}\subset U_{g_1g_2}$ and since $\L$ is $G$-equivariant, we have 
\beq\label{restriction}
(c'^*_{g_1g_2}\L)|_{(U_{g_1}\cap U_{g_2})^\circ}\cong (c'^*_{g_1}\L)|_{(U_{g_1}\cap U_{g_2})^\circ}\otimes (c'^*_{g_2}\L)|_{(U_{g_1}\cap U_{g_2})^\circ}.
\eeq

 Let us recall \cite[Prop. 4.2.]{De}:
\bprop\label{support}
The objects of $e\D_U(G)$ are supported on those $g\in G$ such that ${c'_g}^*\L$ is the trivial  multiplicative local system on $U^\circ_g$.
\eprop

For example, we see that the objects of $e\D_U(U)$ can only be supported on $K\subset U$ where as defined above $K/H=\ker(\phi_\L)=K_\L\subset U/H$.
For $g\in G$, the full subcategory $e\D_U(Ug)\subset e\D_U(G)$ consisting of the complexes supported on the coset $Ug$ is a $e\D_U(U)$ bimodule category. For any $M\in e\D_U(G)$, $M|_{Ug}\in e\D_U(Ug)$. Let $\tM^{perv}_{Ug,e}\subset e\D_U(Ug)$ be the full subcategory consisting of objects whose underlying complex is a perverse sheaf and let $\tM_{Ug,e}=\tM_{\g,e}:=\tM^{perv}_{Ug,e}[\dim H]\subset e\D_U(Ug)$ where $Ug=\g\in G/U=\Gamma$. The following results have been proved in \cite{De}.

\bthm\label{DeMain1}
Using the notation above, we have
\bit
\item[(i)] $\tMg$ is  a semisimple abelian category with finitely many simple objects and $e\D_U(Ug)\cong D^b(\tMg)$. 
\item[(ii)] The simple objects of $\tM_{U,e}$ are precisely the translates $e^k$ of the idempotent $e$ by $k\in K\subset U.$ The isomorphism classes of simple objects correspond to the $H$-cosets in $K$, i.e. to $K/H=\ker(\phi_\L)=K_\L$. In fact $\tM_{U,e}$ has the structure of a modular category equivalent to the modular category $\M(K_\L,\theta)$ associated with the metric group $(K_\L,\theta)$.
\item[(iii)] These abelian categories are well behaved under convolution with compact support, namely we have $\ast: \tMga\times \tMgb\to \tMgab$ for $g_1,g_2\in G$. In particular, each $\tM_{Ug,e}$ is a $\tM_{U,e}$-module category.
\item[(iv)] There is a monoidal action of $\Gamma$ (considered as a discrete group for now) on the monoidal category $e\D_U(G)$ such that $\gamma(\tM_{\g',e})\subset \tM_{\g\g'\g}.$
\item[(v)] For $\g=Ug\in \Gamma, M\in e\D_U(Ug), N\in e\D_U(G)$, we have functorial isomorphisms
\beq
\beta_{M,N}: M\ast N\to \g(N)\ast M
\eeq
which  satisfy certain compatibilities.
\item[(vi)] If $M\in e\D_U(Ug)\subset e\D_U(G)$, then $M$ has a rigid dual $M^\vee\in e\D_U(Ug^{-1})$ and hence we have a natural isomorphisms $\Hom(M^\vee\ast A, B)\cong \Hom(A, M\ast B)$ for $A,B\in e\D_U(G)$.

\eit
\ethm

\subsection{Central Heisenberg admissible pairs}\label{centralheispair}
In this section, let us suppose that the Heisenberg admissible pair $(H,\L)$ is central, i.e. the induced action of $G$ and hence of $\Gamma$ on $U/H$ is trivial.  In this case, we have the commutator map $c:U\times G\to H$, and hence for each $g\in G$, $U_g$ as defined in \S\ref{crossedbraided} is all of $U$. Hence by Proposition $\ref{support}$ in this case we see that the objects of $e\D_U(G)$ are supported on $\tK=\{g\in G|c_g^*\L\cong\Qlcl\}$ which is a closed subgroup $H\subset \tK\subset G$. In fact we have
\bprop\label{centralcrossedbraided} (See \cite[Prop. 2.20]{De13}.)
For each $k\in \tK$, the translates $e^k$ lie in $\tM_{Uk,e}$ and the isomorphism class of $e^k$ only depends on the coset $Hk.$ This classifies all the simple objects of the different semisimple abelian categories $\tMg\subset e\D_U(Ug)$ of objects supported on single cosets.
\eprop

In this situation we get a central extension 
\beq\label{crbrcentral}
0\to K_\L\to K_\Gamma:=\tK/H\to \Gamma\to 0.
\eeq
We have $K=\tK\cap U$ and an extension $0\to K\to \tK\to \Gamma\to 0$.

\section{Idempotents on tori}\label{torus}

In this section, we prove that there are no non-trivial idempotents on a torus:

\bprop\label{main}
Let $e\in \D(T)$ be an idempotent. Then either $e\cong \delta_1$ or $e\cong 0$.
\eprop

This proposition is a consequence of the theory of perverse sheaves on tori developed by O. Gabber, F. Loeser and N. Katz in \cite{GL} and \cite{K}. In this section, we recall this theory and prove the proposition.



Since pushforwards with compact support along group homomorphisms are monoidal functors, they preserve idempotents and we have the following:

\blem\label{cohomology}
Let $\pi:T\to \{1\}$ be the structure morphism and let $e\in \D(T)$ be an idempotent. Then either $\pi_!e\cong \Qlcl$ or $\pi_!e\cong 0$.
\elem

Let $\pi_1^t(T)$ denote the tame fundamental group of the torus and $\pi_1^l(T)$ the maximal pro-$l$ quotient of $\pi_1^t(T)$. Let $\hZ(1):=\lim\limits_{\stackrel{\longleftarrow}{(p,N)=1}}\mu_N(\k)$ and let $\Z_l(1):=\lim\limits_{{\longleftarrow}}\mu_{l^N}(\k)$.
Let $X_*(T):=\Hom(\G_m,T)$ denote the free abelian group of cocharacters of the torus $T$. Let $\hX(T):=X_*(T)\otimes_\Z \hZ(1)$ and $X_*^l(T):=X_*(T)\otimes_\Z \Z_l(1)$. Then we have identifications $\pi^t_1(T)\cong \hX(T)$ and $\pi^l_1(T)\cong X^l_*(T)$.

For a continuous character $\chi:\pi_1^t(T)\cong\hX(T)\to\Qlcl^*$, let $\L_\chi$ denote the corresponding rank 1 Kummer local system on $T$. It is a multiplicative local system, namely we have an isomorphism $\mu^*\L_\chi\cong \L_\chi\boxtimes \L_\chi$. Let us denote the group of continuous characters $\chi:\pi_1^t(T)\cong\hX(T)\to\Qlcl^*$ by $\C(T,\Qlcl^*)$.
Using the fact that $\L_\chi$ defined above is multiplicative, we have:
\blem\label{Ltwist}
Let $\chi\in \C(T,\Qlcl^*)$. Then the functor $\cdot\otimes \L_\chi: \D(T)\to \D(T)$ is monoidal. In particular, if $e\in \D(T)$ is an idempotent, so is $e\otimes \L_\chi$.
\elem

Let $\C^f(T,\Qlcl^*)$ be the subgroup of characters of a finite order coprime to $l$ and let $\C^l(T,\Qlcl^*)$ be the group of continuous characters $\chi:\pi^l_1(T)\to \Qlcl^*$. We have $\C(T,\Qlcl^*)=\C^f(T,\Qlcl^*)\times \C^l(T,\Qlcl^*)$.

For a $\k$-scheme $X$ and an $A\in \D(X):=D^b_c(X,\Qlcl)$, we have an equality of the usual Euler characteristic
$$\chi(X,A):=\sum{(-1)^i\dim H^i(X,A)}$$ and the Euler characteristic with compact support
$$\chi_c(X,A):=\sum{(-1)^i\dim H_c^i(X,A)}.$$ We usually denote the Euler characteristic as just $\chi(X,A)$.

We will use the following result of Deligne:
\bprop\label{Deligne}
For any $A\in\D(T)$ and any $\chi\in\C(T,\Qlcl^*)$ we have $\chi(T,A\otimes\L_\chi)=\chi(T,A)$.
\eprop

\subsection{The Mellin transform}
In \cite{GL}, Gabber and Loeser define a Mellin transform (with compact support) functor $\M_!:\D(T)\to D^b_{coh}(\C(T))$, where $\C(T)$ is a $\Qlcl$-scheme parametrizing the continuous characters of the tame fundamental group of $T$. We briefly recall this below and refer to {\it op. cit.} for more details. Each connected component of $\C(T)$ is isomorphic as a scheme to $\C^l(T)$, the $\Qlcl$-scheme parametrizing the continuous characters of $\pi^l_1(T)$. More precisely,
\beq
\C(T)=\bigsqcup\limits_{\chi\in\C^f(T,\Qlcl^*)}\{\chi\}\times \C^l(T).
\eeq
The closed points, or equivalently the $\Qlcl$-points of $\C^l(T)$ (resp. $\C(T)$) can be identified with the character group $\C^l(T,\Qlcl^*)$ (resp. $\C(T,\Qlcl^*)$).

\brk
When dealing with functors between triangulated categories (of both coherent and $\Qlcl$-constructible type) associated with schemes we often skip writing `L' and `R' and all functors are considered in the derived sense. In particular, inverse image functors between derived categories of coherent sheaves are always the total derived inverse image functors. Similarly all tensor products are derived tensor products.
\erk

\bprop\label{mellin}
As before, let $\pi$ be the structure morphism of $T$. For $A\in \D(T)$ and an inclusion of a closed point $\{\chi\}\stackrel{i_\chi}{\hookrightarrow}\C(T)$ we have a canonical isomorphism $$i_\chi^*\M_!(A)\cong \pi_!(A\otimes \L_\chi).$$
\eprop
In other words, the Mellin transform (with compact support) is a gadget which stores information about the cohomology (with compact supports) of all the Kummer twists of an object of $\D(T)$.

The Mellin transform takes convolution (with compact support) in $\D(T)$ to tensor product in $D^b_{coh}(\C(T))$.
\bprop \cite{GL}.
For $A,B\in \D(T)$, $\M_!(A\ast B)\cong \M_!A\otimes \M_!B$.
\eprop

The following important result is due to Laumon:
\bprop\label{Laumon} (\cite[Prop. 3.4.5]{GL})
Let $\pi:T\to \{1\}$ be the structure morphism. Let $A\in \D(T)$. Then $A=0$ if and only if $\pi_!(A\otimes \L_\chi)=0$ for all $\chi\in \C(T,\Qlcl^*)$. In other words, $A=0$ if and only if $\M_!(A)=0$.
\eprop

As a corollary we can prove that $\D(T)$ has no nilpotents.
\bcor\label{nonil}
Suppose $A\in \D(T)$ is such that $A^{\ast n}:=A\ast\cdots \ast A=0$. Then $A=0$.
\ecor
\bpf
If $A$ is a nilpotent object in $\D(T)$ as above, then $\pi_!A$ is nilpotent in $D^b(Vec)$, hence $\pi_!A=0$. By Lemma \ref{Ltwist} for each $\chi\in \C(T,\Qlcl^*)$, $A\otimes \L_\chi$ is nilpotent and hence $\pi_!(A\otimes \L_\chi)=0$.  Hence by Proposition \ref{Laumon}, A=0.
\epf

\subsection{Euler characteristic zero case}

As a consequence of the results stated above, we can now prove:
\bprop
Let $e\in\D(T)$ be an idempotent. If $\pi_!e=0$, then $e=0$.
\eprop
\bpf
Since $\pi_!e=0$, $\chi(T,e)=0$. Hence by Proposition \ref{Deligne}, $\chi(T,e\otimes\L_\chi)=0$ for all $\chi\in \C(T,\Qlcl^*).$ But by Lemma \ref{Ltwist}, each Kummer twist $e\otimes \L_\chi$ is also an idempotent of Euler characteristic zero. Hence by Lemma \ref{cohomology}, we must have that $\pi_!(e\otimes \L_\chi)=0$ for every $\chi\in \C(T,\Qlcl^*).$ Now using Proposition \ref{Laumon} we see that $e=0$.
\epf

\subsection{Euler characteristic 1 case}
We are now reduced to idempotents of Euler characteristic 1. Let $e$ be an idempotent of Euler characteristic 1. In fact as we have already seen, we must have $\pi_!e\cong \Qlcl$. By Lemmas \ref{cohomology}, \ref{Ltwist} and Proposition \ref{Deligne} we see that we have the following: 
\blem\label{cohomologyoftwists}
If $e\in\D(T)$ is an idempotent of Euler characteristic 1, then $\pi_!(e\otimes \L_\chi)\cong \Qlcl$ for every $\chi\in \C(T,\Qlcl^*)$.
\elem
\bpf
We only need to observe that $e\otimes\L_\chi$ is also an idempotent of characteristic 1.
\epf

\blem
The Mellin transform $\M_!e$ is isomorphic to the free sheaf of rank 1. The idempotent $e$ is in fact a perverse sheaf.
\elem
\bpf
Since Mellin transform takes convolution to tensor product, $\M_!e$ must be an idempotent under tensor product. Using Lemma \ref{cohomologyoftwists}, Proposition \ref{mellin} and the fact that $\M_!e$ is an idempotent, we see that $\M_!e$ must be the free rank 1 sheaf on $\C(T)$. Now using \cite[Thm. 3.4.7]{GL} we conclude that $e$ must be a perverse sheaf.
\epf

Let $\Perv(T)$ be the category of perverse sheaves on $T$. Let $S(T)\subset \Perv(T)$ be the full subcategory consisting of perverse sheaves of Euler characteristic zero. It is a thick subcategory. The following result from \cite{GL} describes all the irreducible perverse sheaves in $S(T)$:
\bprop\label{negligible}(\cite[Thm. 0.3]{GL})
Let $A$ be an irreducible perverse sheaf on $T$ of Euler characteristic zero. Then there exists an exact sequence $0\to\G_m\to T\stackrel{p}{\rightarrow} T'\to 0$, a perverse sheaf $A'$ on $T'$ and a closed point $\chi\in \C(T)$ such that $$A\cong\L_\chi\otimes p^*(A'[1]).$$
\eprop

\bdefn(\cite[Def. 3.7.1]{GL})
For $A\in \Perv(T)$, let $A_t$ be the largest subobject of $A$ lying in $S(T)$ and let $A^t$ be the  subobject of $A$ such that $A/A^t$ is the largest quotient of $A$ lying in $S(T)$. Let $\Perv_{{int}}(T)\subset \Perv(T)$ be the full subcategory formed by perverse sheaves $A$ such that $A_t=0$ and $A/A^t=0$. For any $A\in \Perv(T)$, let $A_{int}:=A^t/(A^t\cap A_t)\cong (A^t+A_t)/A_t\in\Perv_{int}(T)$. This functor $\Perv(T)\to \Perv_{int}(T)$ induces an equivalence $\Perv(T)/S(T)\cong \Perv_{int}(T).$
\edefn

\bprop (\cite[Thm. 3.7.5]{GL})
$\Perv(T)/S(T)$ is a neutral Tannakian category with unit object $\delta_1\in\Perv_{int}(T)\cong \Perv(T)/S(T)$. Hence if $e\in \Perv(T)$ is an idempotent of Euler characteristic 1, then $e_{int}\cong \delta_1$.
\eprop


We will now use induction on the dimension of $T$ to prove that:
\bprop
Let $e\in\D(T)$ be an idempotent such that $\chi(T,e)=1$. Then $e\cong \delta_1$.
\eprop
\bpf
The proof is based on an idea communicated by Antonio Rojas Le\'{o}n. 
The assertion is clear in case $\dim(T)=0$. Let us now assume that this holds for all tori of smaller dimension that $T$. We have already seen that $e$ must be a perverse sheaf and that $e_{int}=\delta_1$. Consider any exact sequence $0\to\G_m\to T\stackrel{p}{\rightarrow} T'\to 0$ and any $A'\in\Perv(T')$. Now $p_!e$ is an idempotent of Euler characteristic 1 on the torus $T'$. Hence by the induction hypothesis, $p_!e$ must be isomorphic to $\delta_1\in \D(T')$. 

We first claim that the perverse sheaf $p^*A'[1]$ (which lies in $S(T)$) does not occur as a quotient of $e$. If possible, let there exist an exact sequence $0\to A\to e\to p^*A'[1]\to 0$ in $\Perv(T)$. Pushing forward with compact support along $p$, we get a distinguished triangle $$\to p_!A\to \delta_1\to A'[1]\otimes p_!\Qlcl\to.$$
Now $p_!\Qlcl$ lives in cohomological degrees 1 and 2 and is a rank one local system in these degrees. The perversity amplitude of the functor $p_!$ is $[0,1]$. Hence taking the assosiated long exact sequence of perverse cohomology (namely the cohomology associated to the perverse $t$-structure), we arrive at a contradiction. 

This shows that for any closed point $\chi\in \C(T)$, $p^*A'[1]$ cannot be the quotient of the idempotent $e\otimes \L_{\chi^{-1}}\in \D(T)$, or in other words that $\L_\chi\otimes p^*(A'[1])$ cannot be a quotient of $e$. By Proposition \ref{negligible}, no nonzero quotient of $e$ can lie in $S(T).$ Hence by \cite[\S2.6]{K}, convolution by $e$ preserves the subcategory $\Perv(T)\subset \D(T)$ and defines an exact functor on $\Perv(T)$.

We have already seen that $e_{int}\cong\delta_1$. Since $e$ has no nonzero quotients in $S(T)$, $e_{int}=e/e_t$ and hence  we must have a surjection $e\twoheadrightarrow \delta_1$ in $\Perv(T)$.   Since convolution with $e$ is exact, we get a surjection $e\onto e$ in $\Perv(T)$ which must be an isomorphism since all objects of $\Perv(T)$ have finite length. This means that $e\onto\delta_1$ is an open idempotent. Then completing to a cone $\to e\to \delta_1\to e'\to$, we get a closed idempotent $\delta_1\to e'$. We see that $e'$ is an idempotent such $\chi(T,e')=0$, hence $e'=0$. Hence we conclude that $e\cong \delta_1$.
\epf

\section{Solvable groups}\label{solvable}
In this section we prove Theorems \ref{main1}, \ref{main2}  and  \ref{main3}.
\subsection{Heisenberg idempotents on connected solvable groups}\label{heisidsol}
Let us first study Heisenberg idempotents on connected solvable groups. We will use notation defined in \S\ref{anheisid}. Let $G$ be a connected solvable group with Heisenberg admissible pair $(H,\L)$. As before, let $U$ be the unipotent radical and since $G$ is solvable, $G/U=:T$ is a torus. Moreover, we have a (non-unique) embedding $T\subset G$ which allows us to write $G=UT.$ All maximal tori in $G$ are conjugate to $T$. The maximal torus $T\subset G$ acts on $U$ and the Heisenberg admissible pair $(H,\L)$ is stabilized by this action. Let $V=U/H$. We have an action of $T$ on $V$ and the $T$-equivariant skewsymmetric isogeny $\phi_\L:V\to V^*$ with associated metric group $(K_\L,\theta)$. By Proposition \ref{skewsymm}, we get a direct sum decomposition of this into:
\bit
\item $\phi_\L^T:V^T\to (V^T)^*$ which has kernel $K_\L$ and
\item $\phi_\L':V'\stackrel{\cong}{\rightarrow}(V')^*.$
\eit

Let $U_T, U' \subset U$ be such that $U_T/H=V^T$  and $U'/H=V'$. It is clear that $(H,\L)$ is an Heisenberg admissible pair for $U_T$ as well as $U'$ with the corresponding isogenies being $\phi_\L^T$ and $\phi'_\L$ respectively.

Note that for $t\in T$, the subgroup $V^t\subset V$ of $t$-fixed points is connected since $t$ is semisimple and $V$ is unipotent (see \cite[\S18]{H}). Hence the subgroup $U_t\subset U$ (see \S\ref{anheisid}) is also connected and we have the map $c'_t:U_t\to H$. As stated before, we may consider $c'^*_t\L$ as a $C_G(t)$-equivariant multiplicative local system on $V^t$.

\blem\label{trivialrestriction}
For $t\in T$, $(c'^*_t\L)|_{V^T}\cong \Qlcl$.
\elem
\bpf
By (\ref{restriction}), for $t_1,t_2\in T$ we have $(c'^*_{t_1t_2}\L)|_{V^T}\cong (c'^*_{t_1}\L)|_{V^T}\otimes (c'^*_{t_2}\L)|_{V^T}$. Hence $T\ni t\mapsto (c'^*_t\L)|_{V^T}$ defines a morphism of perfect group schemes $T\to (V^T)^*$. Since $T$ is a torus and $(V^T)^*$ is unipotent, this map must be trivial. Hence $(c'^*_t\L)|_{V^T}\cong \Qlcl$ for each $t\in T$.
\epf

\bprop\label{TintK}
For $t\in T$, the multiplicative local system $c'^*_t\L$ on $V^t$ is trivial.
\eprop
\bpf
Since $T\subset C_G(t)$, the multiplicative local system $c'^*_t\L$ is $T$-equivariant. By Lemma \ref{fp}, it suffices to prove that the restriction of the multiplicative local system to $(V^t)^T=V^T$ is trivial. We have already proved this in Lemma \ref{trivialrestriction}.
\epf

\subsection{Extension of Heisenberg admissible pairs on connected solvable groups}\label{extonsol}

Now suppose $H\subset L\subset U$ is such that $L/H\subset U/H$ is $T$-stable isotropic for the skew-symmetric biextension $\phi_\L$. Then as in \S\ref{eoap} we can extend the Heisenberg admissible pair $(H,\L)$ to an admissible pair $(L,\L')$. By (\ref{normalizer}), we have an extension
\beq\label{tnormalizer}
0\to L^\perp\to \N_G(L,\L')\to T\to 0.
\eeq
Since the neutral connected component $\N^\circ_G(L,\L')$ is a connected solvable group, the set of its unipotent elements is a connected subgroup. Hence we must have $\N^\circ_G(L,\L')\cap L^\perp=L^{\perp\circ}$ and an embedding $T\hookrightarrow \N^\circ_G(L,\L')\subset \N_G(L,\L')\subset G$ giving splittings $\N^\circ_G(L,\L')\cong L^{\perp\circ}\rtimes T$ and $\N_G(L,\L')\cong L^{\perp}\rtimes T$. In particular, we have an embedding $T\subset G$ which normalizes the pair $(L,\L')$. Since all embeddings of $T$ in $G$ are conjugate, there must exist $u\in U\subset G$ such that the torus ${ }^uT$ normalizes the pair $(L,\L')$ where $T\subset G$ is the chosen embedding that we started out with. This implies that $T$ normalizes the admissible pair $(L,{ }^{u^{-1}}\L')$. Hence we have proved:
\blem\label{solext}
Let $G=UT$ be a connected solvable group with a chosen maximal torus $T\subset G$ and a Heisenberg admissible pair $(H,\L)$ for $G$. Suppose $H\subset L\subset U$ is such that $L/H\subset U/H$ is a $T$-stable isotropic subgroup for the skew-symmetric biextension $\phi_\L$. Then there exists an  admissible pair $(L,\L')$ that extends $(H,\L)$ and is normalized by $T\subset G$. In this case $\N_G(L,\L')= L^{\perp} \cdot T$.
\elem

\subsection{Central admissible pairs for connected solvable groups}\label{centralsol}
Now suppose that the Heisenberg admissible pair $(H,\L)$ for the connected solvable group $G$ is central. We use all the notation from \S\ref{centralheispair}. By Proposition \ref{TintK}, we see that $T\subset \tK$ i.e. $c_t^*\L\cong\Qlcl$ and we have splittings $\tK=KT$ and $K_T\cong K_\L\times T$ (see (\ref{crbrcentral})). Since the action of $T$ on $U/H$ is trivial, we see that $G/H$ splits as a direct product $G/H\cong U/H\times T$. Note that using the isomorphism $c_t^*\L\cong\Qlcl$ of multiplicative local systems, for $t\in T, M\in e\D_U(G)$ we have functorial isomorphisms $M\xto{\cong}t(M)$ which in fact give a trivialization of the action of $T$ on $e\D_U(G)$. 
We now prove a special case of Theorem \ref{main1}:
\bprop\label{specialmain1}
In this situation we have a canonical equivalence of monoidal categories $$ e\D_U(U)\boxtimes \D(T)\xto{\cong}e\D_U(G)$$ which is independent of the choice of the splitting $T=G/U\hookrightarrow G$.
\eprop
\bpf
We have the chosen embedding $T\subset G$ using which we have $U\times T\cong G$ (via the multiplication map) as a quasi-algebraic scheme. Using this identification, for $M\in \D(U)\subset \D(G)$ and $N\in \D(T)\subset \D(G)$, we have natural isomorphism $M\boxtimes N\cong M\ast N\in \D(G)$. Now suppose that $M\in e\D_U(U)$ and $N\in \D(T)$, then it is clear that $M\ast N\in e\D(G)$ and also easy to check that $M\boxtimes N\in \D_U(G)$. Hence we see that we have a well defined tautological functor $ e\D_U(U)\boxtimes \D(T)\to e\D_U(G)\cong e\D_G(G).$ The multiplication map $U\times T \times U\times T\to U\times T$ is given by $(u_1,t_1,u_2,t_2)\mapsto (u_1\cdot { }^{t_1}u_2, t_1t_2)$. Since $T$ acts trivially on $e\D_U(U)$, we can see that this tautological functor has a monoidal structure. We have seen that the objects of $e\D_U(G)$ have support in $\tK$. Note that $\tK^\circ=K^\circ T=HT$ and the connected components of $\tK$ are the $HT$-cosets $HTk$ in $\tK$. We have an equivalence $e\D(HTk)= \D_{H,\L}(HTk)\cong \D(T)$ of triangulated categories.  Hence we see that this tautological functor is in fact an equivalence of monoidal categories.

Next we check that this is independent of the choice of the embedding $T\subset G$. First note that in the connected solvable group $G$, the normalizer of any maximal torus $T$ is exactly the same as its centralizer. Indeed, it suffices to prove that $N_U(T)=C_U(T)$. Now suppose that $u\in N_U(T)$, i.e. $utu^{-1}=u\cdot tu^{-1}t^{-1}\cdot t\in T$ for each $t\in T$. Since $u\cdot tu^{-1}t^{-1}\in U$, we must have $u\cdot tu^{-1}t^{-1}=1$ for each $t\in T$, i.e. $u\in C_U(T)$. 

Next, recall that we have $G/H\cong U/H\times T$ and hence $G/H$ has a unique maximal torus. The preimage of this maximal torus in $G$ is equal to $HT\normal G$. Hence all maximal tori in $G$ are conjugate under the action of $H\subset G$. Moreover we have noted that the normalizer of a torus equals it centralizer. Hence we conclude that any splitting $T=G/U\hookrightarrow G$ can be obtained from our chosen embedding $T\subset G$ by conjugation by some element $h\in H\subset G$. Using the maximal torus ${}^{h}T\subset G$, we obtain the new identification $G\cong U\times {}^hT$ as varieties. The induced isomorphism $U\times T\rar\cong U\times {}^hT$ is defined by $(u,t)\mapsto (u\cdot tht^{-1}\cdot h^{-1},{}^ht)$. We see that for $M\in e\D_U(U)$ and $N\in \D(T)$, the pullpack of the object $M\boxtimes {}^hN$ on $U\times {}^hT$ by the above isomorphism is functorially isomorphic to the object $M\boxtimes N$ on $U\times T$. This constructs a natural isomorphism between the monoidal equivalence $e\D_U(U)\boxtimes \D(T)\xto{\cong}e\D_U(G)$ obtained from our chosen maximal torus $T\subset G$ and the composition $e\D_U(U)\boxtimes \D(T)\xto{\id\boxtimes {}^h(\cdot)} e\D_U(U)\boxtimes \D({}^hT)\xto{\cong}e\D_U(G)$ which corresponds to the splitting $T\xto{\h{ad}(h)}{}^hT\subset G$.
\epf

\brk
The category $e\D_U(U)\boxtimes\D(T)$ has a (trivial) braided $T$-crossed structure and the above equivalence is compatible with this extra structure. 
\erk

\subsection{Proof of Theorem \ref{main1}}\label{pfmain1}
We continue to use the notation of \S\ref{heisidsol}. As before, we have a Heisenberg admissible pair $(H,\L)$ on $G=UT$. In this section, we no longer assume that the action of $T$ on $V=U/H$ is trivial. As we have seen in \S\ref{heisidsol}, we have a canonical splitting $V=V^T\oplus V'$.  By Theorem \ref{stabisot}, there exists a $T$-stable Lagrangian subgroup in $V'$. (It will be isotropic in $V$.) Let $H\subset L\subset U' \subset U$ be such that $L/H\subset V'$ is $T$-stable Lagrangian subgroup in $V'$. We have $(L/H)^\perp=V^T\oplus (L/H)\subset V$. Then by Lemma \ref{solext}, we have an admissible pair $(L,\L')$ normalized by $T$ such that $(H,\L)\leq (L,\L')$. We have $\N_G(L,\L')=U_TLT$, $\R_u(\N_G(L,\L'))=U_TL$ and $U_TL/L=U_T/H=V^T$. Hence $(L,\L')$ is a central admissible pair. By Proposition \ref{specialmain1}
\beq
e_{\L'}\D_{U_TL}(U_TL)\boxtimes \D(T)\xto{\cong}e_{\L'}\D_{U_TL}(U_TLT).
\eeq
By Proposition \ref{heisext} (and its method of proof) we get equivalences 
\beq 
\ind_{U_TL}^U: e_{\L'}\D_{U_TL}(U_TL)\xto{\cong} e\D_U(U) \hbox{ of modular categories and}
\eeq
\beq\label{crbreq}
\av_{{U}/{U_TL}}: e_{\L'}\D_{U_TL}(U_TLT)\xto{\cong} e\D_U(G) \hbox{ of monoidal categories}
\eeq
with a quasi-inverse given by $M\mapsto e_{\L'}\ast M$. 
\brk
The equivalence (\ref{crbreq}) is in fact an equivalence of braided $T$-crossed categories.
\erk
Thus we get an equivalence \beq\label{crbrequiv} e\D_U(U)\boxtimes \D(T)\xto{\cong} e_{\L'}\D_{U_TL}(U_TL)\boxtimes \D(T)\xto{\cong}e_{\L'}\D_{U_TL}(U_TLT)\xto\cong e\D_U(G)\eeq
given by $M\boxtimes N\mapsto \av_{U/U_TL}((e_{\L'}\ast M)\boxtimes N)\xto\cong \av_{U/U_TL}(e_{\L'}\ast M\ast N)$ for $M\in e\D_U(U)$ and $N\in \D(T).$ 

Note that we have used the choice of the central admissible pair $(L,\L')$ to define this monoidal equivalence. Next we will prove that the monoidal equivalence does not depend on this choice.

Suppose that $(L_1,\L_1)$ and $(L_2,\L_2)$ are two such pairs. In particular $L_1/H$ and $L_2/H$ are $T$-stable Lagrangian subgroups in $V'$ and $L_1,\L_2$ are $T$-equivariant multiplicative local systems. Let $e_i=\L_i\otimes \can_{L_i}\in \D_{U_TL_i}(U_TL_iT)$ be the corresponding Heisenberg idempotent. Note that since $V'^T=0$, we must have $\left(\frac{L_1\cap L_2}{H}\right)^T=0$ and hence ${\left(\frac{L_1\cap L_2}{H}\right)^*} ^T=0$ by Lemma \ref{l:teqls}. This implies that $\L_1|_{L_1\cap L_2}\cong \L_2|_{L_1\cap L_2}$. Let $\L'$ denote this restriction to the intersection. Then $(L_1\cap L_2,\L')$ is also an admissible pair for $G$ by \S\ref{eoap} with normalizer $U_TL_1L_2T\subset G$. Let $e'=\L'\otimes \can_{L_1\cap L_2}$ be the corresponding idempotent.

Now for $i=1,2$, the equivalence $\av_{U/U_TL_i}: e_i\D_{U_TL_i}(U_TL_iT)\rar\cong e\D_U(G)$ is naturally isomorphic to the composition $e_i\D_{U_TL_i}(U_TL_iT)\xto{\av{U_TL_1L_2/U_TL_i}} e'\D_{U_TL_1L_2}(U_TL_1L_2T)\xto{\av_{U/U_TL_1L_2}} e\D_U(G)$. Now we want to show that the compositions $\D(T)\rar{}e_i\D_{U_TL_i}(U_TL_iT)\rar\cong e\D_U(G)$ for $i=1,2$ are naturally isomorphic. For this it is enough to show that for $i=1,2$, the compositions $\D(T)\rar{}e_i\D_{U_TL_i}(U_TL_iT)\rar\cong e\D_{U_TL_1L_2}(U_TL_1L_2T)$ are naturally isomorphic.

Essentially we have reduced the problem to the situation where $L_1/H$ and $L_2/H$ are complementary $T$-stable Lagrangian subgroups of $V'$ (by passing from $G$ to be the subgroup $U_TL_1L_2T$ and from the pair $(H,\L)$ to the pair $(L_1\cap L_2,\L')$ from the previous paragraph). So let us assume without loss of generality that $L_1\cap L_2=H$ and hence that $U'=L_1L_2$ (see \S\ref{heisidsol}).

Now it suffices to prove that the triangulated monoidal functor from $\D(T)$ to $e_2\D_{U_TL_2}(U_TL_2T)$ defined by (cf. Proposition \ref{weakid})
$$X\mapsto e_2\ast\av_{U/U_TL_1}(e_1\boxtimes X)$$ is naturally isomorphic to the canonical triangulated monoidal functor $\D(T)\rar{}e_2\D_{U_TL_2}(U_TL_2T)$ defined by $X\mapsto e_2\boxtimes X$ (cf. Proposition \ref{specialmain1}). Now $U_TL_1$ is a normal subgroup of $U$ and we have identifications of the quotients $U/U_TL_1=U'/L_1=L_2/H$. Moreover we have an action of this quotient on $\D_{U_TL_1}(U_TL_1T)$ induced by conjugation. To compute $e_2\ast\av_{U/U_TL_1}(e_1\boxtimes X)$ we find it convenient to use the `integral' notation ($\int$) for pushforward with compact support.

Now $e_2\ast\av_{U/U_TL_1}(e_1\boxtimes X)$ is an object of $e_2\D_{U_TL_2}(U_TL_2T)$ and hence is only supported on $L_2T\subset G=UT$. For $t\in T, v\in L_2\subset U$, the stalk of $(e_2\ast\av_{U/U_TL_1}(e_1\boxtimes X))$ at $vt$ is given by
$$(e_2\ast\av_{U/U_TL_1}(e_1\boxtimes X))(vt)=\int\limits_{l_2\in L_2}e_2(l_2^{-1})\otimes\av_{U/U_TL_1}(e_1\boxtimes X)(l_2vt)$$
$$=\int\limits_{l_2\in L_2}\int\limits_{uH\in L_2/H}e_2(l_2^{-1})\otimes(e_1\boxtimes X)(ul_2vtu^{-1})$$
$$=\int\limits_{l_2\in L_2}\int\limits_{uH\in L_2/H}e_2(l_2^{-1})\otimes e_1(ul_2vtu^{-1}t^{-1})\otimes X(t)$$
$$=X(t)\otimes\int\limits_{l_2\in L_2}\int\limits_{uH\in L_2/H}e_2(l_2^{-1})\otimes e_1(ul_2vtu^{-1}t^{-1}).$$
For any $u\in L_2$, let us set $\bar{u}=uH\in L_2/H=U/L_1$. Now the object $e_1$ is only supported on $L_1$. So we see that the above stalk is canonically isomorphic to 
$$X(t)\otimes\int\limits_{\substack{l_2\in L_2, \bar{u}\in L_2/H\\(t-\id)\bar{u}=\bar{l_2v}}}e_2(l_2^{-1})\otimes e_1(ul_2vtu^{-1}t^{-1})$$
$$=X(t)\otimes\int\limits_{\substack{\bar{u}\in L_2/H, h\in H\\{(l_2v=htut^{-1}u^{-1})}}}e_2(vutu^{-1}t^{-1}h^{-1})\otimes e_1(u\cdot htut^{-1}u^{-1}\cdot tu^{-1}t^{-1})$$
$$=X(t)\otimes\int\limits_{{\bar{u}\in L_2/H, h\in H}}e_2(v)\otimes \L_2(u)\otimes \L_2(tu^{-1}t^{-1})\otimes \L_2(h^{-1}))\otimes \L_1(uhu^{-1})\otimes \L_1([u,tut^{-1}])\otimes e_1(1).$$
Now $\L_2$ is a $T$-equivariant multiplicative local system on $L_2$. Moreover, $L_2/H\subset U/H$ is a $T$-stable Lagrangian subgroup. Hence the pullback of $\L=\L_1|_H$ via the commutator $L_2\times L_2\rar{[,]}H$ is trivial. Hence the above stalk is canonically isomorphic to 
$$e_2(v)\otimes X(t)=(e_2\boxtimes X)(vt)$$
as desired. This completes the proof that the monoidal equivalence $e\D_U(U)\boxtimes \D(T)\cong e\D_U(G)$ that we have defined does not depend on any choices.

\subsection{Groups with solvable neutral connected component}
In this section we study Heisenberg admissible pairs on possibly disconnected groups $G$ whose neutral connected component $G^\circ$ is solvable. Let $(H,\L)$ be the Heisenberg admissible pair on $G$ and let $e\in \D_G(G)$ be the corresponding Heisenberg idempotent. As before let $U:=\R_u(G)=\R_u(G^\circ)$. Let $G^\circ/U=T$. We identify $T$ with a fixed maximal torus in $G^\circ$ (i.e. we fix lifts in $G^\circ$ for each element of $G^\circ/U$ such that these lifts define a maximal torus) using which we have $G^\circ=UT$. As before, let $\Gamma=G/U$. We have $\Gamma^\circ=T$. We have an action of $\Gamma$ (and hence also of $G$) on $T$ by conjugation. Since the action of $T$ is trivial, we get an action of $\Pi_0:=\pi_0(\Gamma)=\pi_0(G)=G/G^\circ$ on $T$. Hence we get a braided monoidal action of $G,\Gamma, \Pi_0$ on the symmetric monoidal category $\D(T)$. 

Let us consider the monoidal Hecke subcategory $e\D_U(G)$ which has the $\Pi_0$-grading $$e\D_U(G)=\bigoplus\limits_{G^\circ g \in G/G^\circ} e\D_U(G^\circ g=UTg)$$
with identity component $e\D_U(G^\circ).$ We will now describe each graded component $e\D_U(G^\circ g)$ as a $e\D_U(G^\circ)$-bimodule category.

We have a canonical equivalence $e\D_U(U)\boxtimes \D(T)\xto{\cong} e\D_U(G^\circ)$ of triangulated monoidal categories given by (\ref{crbrequiv}) which in particular gives a canonical fully faithful monoidal functor $\D(T)\to e\D_U(G^\circ)$ given by $\D(T)\ni N\mapsto \bar{N}:=\av_{U/U_TL}(e_{\L'}\ast N)=\av_{U/U_TL}(e_{\L'}\boxtimes N)\in e\D_U(G^\circ)$ using the notation from (\ref{crbrequiv}), which identifies $\D(T)$ as a direct summand of $e\D_U(G^\circ)$ and we have a left quasi-inverse $e\D_U(G^\circ)\xto{\phi} \D(T)$. Moreover, we have seen that this functor does not depend on the choice of the central admissible pair $(L,\L')$ for $G^\circ$.  Hence we see that the functor $\D(T)\to e\D_U(G^\circ)$ is compatible with the actions of $\Gamma$ on $\D(T)$ and $e\D_U(G^\circ)$.

 By Theorem \ref{DeMain1}, we know that for each $g\in G$, $e\D_U(Ug)\cong D^b(\tM_{Ug,e})$ is an $e\D_U(U)$-bimodule category and that $\tM_{Ug,e}$ is a semisimple abelian category with finitely many simple objects.
\bdefn
For $g\in G$, we define a  $e\D_U(U)\boxtimes \D(T)$-bimodule structure on  $e\D_U(Ug)\boxtimes \D(T)$ by setting 
\beq
(M\boxtimes N)\ast (A\boxtimes K):=(M\ast A)\boxtimes (N\ast K) \in e\D_U(Ug)\boxtimes \D(T)
\eeq 
\beq
(A\boxtimes K)\ast (M\boxtimes N):=(A\ast M)\boxtimes (g^{-1}(K)\ast N) \in e\D_U(Ug)\boxtimes \D(T)
\eeq 
for $M\boxtimes N\in e\D_U(Ug)\boxtimes \D(T)$ and  $A\boxtimes K\in e\D_U(U)\boxtimes \D(T)$. 
\edefn

We will now prove
\bprop
We have an equivalence $e\D_U(Ug)\boxtimes \D(T)\xrightarrow\cong e\D_U(G^\circ g)$ of $e\D_U(G^\circ)$-bimodule categories given by $M\boxtimes N\mapsto M\ast \bar{N}\in e\D_U(G^\circ g)$ for $M\in e\D_U(Ug), N\in \D(T)$. 
\eprop
\bpf
Using Theorem \ref{DeMain1}, we can check that the functor above has a natural $e\D_U(G^\circ)$-bimodule functor structure. Also, we know from Theorem \ref{DeMain1} that $e\D_U(Ug)$ is the bounded derived category of a semisimple abelian category and that every object $M\in e\D_U(Ug)$ has a rigid dual $M^\vee\in e\D_U(Ug^{-1})$. Let $M$ be a simple object in $e\D_U(U_g)$. Then by \cite[Prop. 5.4]{De}, $M^\vee\ast M\cong e\oplus\left(\left(\oplus e^k \right)\right)$, a direct sum of $e$ and distinct translates of the idempotent $e$ and we have a canonical evaluation map $M^\vee\ast M\to e$.  Consider the functor $\D(T)\to e\D_U(G^\circ)\to e\D_U(G^\circ g)$ defined by $N\mapsto M\ast \bar{N}$. This has a left quasi-inverse given by $e\D_U(G^\circ g)\ni X\mapsto \phi(M^\vee \ast X)$, i.e. the composition 
\beq
e\D_U(G^\circ g)\xto{M^\vee\ast(\cdot)} e\D_U(G^\circ) \xto{\phi} \D(T).
\eeq
Using this, we see that we get an equivalence $e\D_U(Ug)\boxtimes \D(T)\xrightarrow\cong e\D_U(G^\circ g)$ of $e\D_U(G^\circ)$-bimodule categories.

\epf

Hence for $g_1,g_2\in G$, we have equivalences 
$$e\D_U(Ug_1)\boxtimes \D(T)\xrightarrow\cong e\D_U(G^\circ g_1),$$
$$e\D_U(Ug_2)\boxtimes \D(T)\xrightarrow\cong e\D_U(G^\circ g_2),$$
$$e\D_U(Ug_1g_2)\boxtimes \D(T)\xrightarrow\cong e\D_U(G^\circ g_1g_2).$$
We see that
\bprop\label{convg1g2}
Under the identifications above, the convolution $e\D_U(G^\circ g_1)\times e\D_U(G^\circ g_2)\to e\D_U(G^\circ g_1g_2)$ is given by the functor 
$$e\D_U(Ug_1)\boxtimes \D(T)\times e\D_U(Ug_2)\boxtimes \D(T)\to e\D_U(Ug_1g_2)\boxtimes \D(T)$$ defined by 
\beq
(M_1\boxtimes N_1)\ast (M_2\boxtimes N_2):=(M_1\ast M_2)\boxtimes (g_2^{-1}(N_1)\ast N_2).
\eeq
\eprop

Note that $T$ acts on each connected component $G^\circ g$ giving rise to an action of $T$ on $e\D_{U}(G^\circ g).$ We will now describe the corresponding action of $T$ on $e\D_U(Ug)\boxtimes \D(T)$ (which is equivalent to $e\D_U(G^\circ g)$.)

Recall that we have the monoidal functor $\D(T)\to e\D_U(G^\circ)$. Under this functor suppose the constant sheaf ${\Qlcl}_T$ on $T$ maps to $E\in e\D_U(G^\circ)$. For $t\in T$, the delta sheaf $\delta_t$ supported at $t$ maps to $E_t:= E|_{Ut}\in e\D_U(Ut)\subset e\D_U(G^\circ)$. Note that for $t_1,t_2\in T$, we have natural isomorphisms $E_{t_1}\ast E_{t_2}\cong E_{t_1t_2}$ and $E_1\cong e$. For $g\in G$, let $\g=Ug\in \Gamma=G/U$. By Theorem \ref{DeMain1}(v), for $M\in e\D_U(Ug)$, we have 
$$e\D_U{\left(Utgt^{-1}\right)}\ni t(M)\cong E_t\ast M \ast E_{t^{-1}}\cong M\ast \g^{-1}(E_t)\ast E_{t^{-1}}$$
$$\cong M\ast E_{\g^{-1}t\g}\ast  E_{t^{-1}}\cong M\ast E_{\g^{-1}t\g t^{-1}}.$$ 

For $\g\in \Gamma$, we have the conjugation automorphism $\g:T\to T$. We have the $\g$-conjugation action $\alpha_\g:T\times T\to T$ defined by $\alpha_\g(t,s)=ts\g(t^{-1})=ts\g t^{-1}{\g^{-1}}$. Thus we get the $\g$-conjugation action of $T$ on $\D(T)$, namely for each $t\in T$ we set the functor $\alpha_\g(t):\D(T)\to \D(T)$ to be the pullback along the map $\alpha_\g(t^{-1},\cdot):T\to T$. For $t\in T$ and $N\in \D(T)$, we have natural isomorphims $\alpha_\g\left(\g^{-1}(t)\right)(N)\cong \delta_{\g^{-1}t\g}\ast N\ast \delta_{t^{-1}}$. Let $\D^\g_T(T)$ denote the $T$-equivariant derived category for the $\g$-conjugation of $T$ on itself.
We see that
\bprop
For $M\in e\D_U(Ug), N\in \D(T)$ and $t\in T$,  set 
$$t(M\boxtimes N):=M\boxtimes \left(\alpha_\g\left(\g^{-1}(t)\right)N\right).$$
With this action of $T$ on $e\D_U(Ug)\boxtimes \D(T)$, the equivalence $e\D_U(Ug)\boxtimes \D(T)\cong e\D_U(G^\circ g)$ respects the $T$-actions. We have an equivalence $e\D_U(Ug)\boxtimes \D^\g_T(T)\cong e\D_{G^\circ}(G^\circ g).$
\eprop

\subsection{Proof of Theorem \ref{main2}}\label{pfmain2}
Suppose $G$ is such that $G^\circ$ is solvable and if $G'$ is a closed subgroup of $G$, then clearly $G'^\circ$ is also solvable. Let $(H,\L)$ be an admissible pair for $G$ with normalizer $G'$. We want to prove that $e_\L\in \D_{G'}(G')$ is a minimal idempotent and $f_\L\in \D_G(G)$ is a minimal idempotent. By Proposition \ref{weakid}, these two assertions are equivalent, so it suffices to prove $e_\L\in \D_{G'}(G')$ is minimal. Without loss of generality assume that $G=G'$, and let $e:=e_\L$. Hence we are now in the situation studied in the previous section. To prove minimality of $e$, we must prove that the Hecke subcategory $e\D_G(G)$ has no idempotents other than $0$ and $e$. For this, it suffices to show that $e\D_{G^\circ}(G)$ has no non-trivial idempotents. 
We first prove that $e\D_{G^\circ}(G)$ has no nilpotents:
\blem\label{nonilsol}
Let $K\in e\D_{G^\circ}(G)$ be such that $K^{\ast n}:=K\ast\cdots \ast K=0$ for some integer $n\geq 1$. Then $K=0.$
\elem
\bpf
Without loss of generality, let $K$ be indecomposable. Hence there is a $g\in G$ such that $K\in e\D_{G^\circ}(G^\circ g)\cong e\D_U(Ug)\boxtimes \D^\g_T(T)$. Since $K$ is indecomposable, $K$ corresponds to some $M\boxtimes N\in e\D_U(Ug)\boxtimes \D^\g_T(T)$. Since $N\in \D^\g_T(T)$, we have natural isomorphisms of stalks $N_s\cong N_{ts\g(t^{-1})}$ for $s,t\in T$ obtained from the $T$-equivariant structure on $N$. Restricting the $T$-equivariant structure isomorphism to the antidiagonal $t=s^{-1}$ we conclude that $\g(N)\cong N$. Hence using Proposition \ref{convg1g2}, we see that $(M\boxtimes N)^{\ast n}\cong (M^{\ast n})\boxtimes (N^{\ast n})=0$. Now since objects of $e\D_U(Ug)$ have rigid duals, $e\D_U(Ug)$ has no nilpotents. Also we have seen that $\D(T)$ has no nilpotents (see Corollary \ref{nonil}). Hence we conclude that either $M=0$ or $N=0$, i.e. $K=0$.
\epf

Finally we prove that 
\bthm
Let $A\in e\D_{G^\circ}(G)$ be an idempotent. Then either $A\cong e$ or $A=0$.
\ethm
\bpf
Let $A=\oplus A_i$ be a decomposition into indecomposables. Since $A_i\ast A_i\neq 0$, we conclude that each $A_i$ must be an idempotent and that $A_i\ast A_j=0$. Now suppose that $A$ is an indecomposable idempotent. Hence we must have $A\in e\D_{G^\circ}(G^\circ)\cong e\D_U(U)\boxtimes \D_T(T)$. Now we have seen in \S\ref{torus} that $\D(T)$ has no non-trivial idempotents and by \cite{BD08}, $e\D_U(U)$ has no non-trivial idempotents. Since $A$ is indecomposable, we must have $A\cong e$. Hence we see that $e\D_{G^\circ}(G)$ has a unique indecomposable idempotent. Hence we conclude that $e, 0$ are the only idempotents in $e\D_{G^\circ}(G)$.
\epf

Hence we conclude that $e\DG$ also has no non-trivial idempotents. Hence the idempotent $e$ obtained using a Heisenberg admissible pair $(H,\L)$ for $G$ is a minimal idempotent in $\DG$ as we had set out to prove.

\subsection{Proofs of Theorem \ref{main3} and Theorem \ref{main4}}\label{pfmain3}
Let $G$ be a group such that $G^\circ$ is solvable. Let $(H,\L)$ be an admissible pair for $G$ with normalizer $G'$. We have proved Theorem \ref{main2}, which shows that $f_\L\in e\DG$ is a minimal idempotent. Thus in this case, we have proved statement (ii) of Conjecture \ref{conjmain}.

To complete the proof of Theorem \ref{main3}, it is sufficient to prove Theorem \ref{main4}, namely that statement (ii) of Conjecture \ref{conjmain} implies statements (i) and (iii).

\bpf[Proof of Theorem \ref{main4}(i)]
Note that by Theorem \ref{compadmipair} there exists an admissible pair $(H,\L)$ for $G$ compatible with $f\neq 0$, i.e. such that $e_\L\ast f\neq 0.$ Using Proposition \ref{weakid} (iv) and (vi), we see that $f_\L\ast f\neq 0.$ By assumption both are minimal idempotents in $\DG$, hence $f\cong f_\L$.
\epf
\bpf[Proof of Theorem \ref{main4}(ii)]
We will prove (ii) by induction on the size of $G$. Let $(H,\L)$ be an admissible pair for $G$. It is sufficient to prove that the idempotent $f_\L$ is locally closed, since by Proposition \ref{weakid}(x) it will then follow that $\indg:e_\L\DGp\to f_\L\DG$ is a braided equivalence. Let $(A,\E)\in \Pn(G)$ be a maximal normal pair compatible with $f_\L$, i.e. $e_\E\ast f_\L\neq 0$. Let $G_1=\N_G(A,\E)$. By Proposition \ref{maxnormpair}, if $G_1=G$, then $(A,\E)$ is a Heisenberg admissible pair for $G$. Hence in this case $e_\E=f_\E\cong f_\L$ must be a closed and in particular a locally closed idempotent.

Now suppose that $G_1$ is a strictly smaller subgroup of $G$ and hence let us assume that the statement holds for $G_1$. By Corollary \ref{loccloscor}, we have a braided equivalence 
\beq\label{equ}
\ind_{G_1}^G: e_\E\D_{G_1}(G_1)\xto{\cong} f_\E\DG
\eeq

and $f_\E$ is a locally closed idempotent in $\DG$. By Proposition \ref{weakid}(iv)
$$f_\E\ast f_\L\cong \ind_{G_1}^G(e_\E\ast f_\L)\neq 0.$$ Since $f_\L$ is a minimal idempotent, we must have $f_\E\ast f_\L\cong f_\L.$ Hence we conclude that $e_\E\ast f_\L$ must be a minimal idempotent in $\D_{G_1}(G_1)$. Using (i) and the induction hypothesis, we conclude that $e_\E\ast f_\L$ is a locally closed idempotent in $\D_{G_1}(G_1)$ and hence in $e_\E\D_{G_1}(G_1)$. By the braided equivalence (\ref{equ}), $f_\L$ must be locally closed in $f_\E\DG$. Since $f_\E$ is locally closed in $\DG$, $f_\L$ must also be locally closed in $\DG$ by Corollary \ref{lcinlc}.
\epf

\section{Towards a theory of character sheaves on algebraic groups}\label{redtoheis}
In this section, we outline a plausible road map towards a theory of character sheaves on arbitrary affine algebraic groups. In particular, we will reduce this problem to the Heisenberg case. We hope to study the Heisenberg case in detail and implement this road map rigorously in future works. Here we only sketch a rough outline.

Let $G$ be any affine algebraic group. Our goal is to define a set $CS(G)$ of character sheaves on $G$ as some special (isomorphism classes of) objects in the braided triangulated category $\DG$. Let $f\in \DG$ be a minimal idempotent. Then we form the Hecke subcategory $f\DG\subset \DG$ and we hope to define the set $CS_f(G)$ of character sheaves on $G$ corresponding to the minimal idempotent $f\in \DG$. To define $CS(G)$ we take the union\footnote{It is easy to see that this will in fact be a disjoint union.} of $CS_f(G)$ over all minimal idempotents $f$ in $\DG$. Now according to Conjecture \ref{conjmain}, any minimal idempotent $f$ should come from some admissible pair $(H,\L)$ for $G$ with normalizer $G'$. Moreover, we expect to have a braided equivalence $\indg:e_\L\DGp\to f\DG$.

Hence to define character sheaves in the piece $f\DG$, we may first define character sheaves in $e_\L\DGp$. Now $(H,\L)$ is a Heisenberg admissible pair for $G'$. Hence to define character sheaves in general, it is enough to study the category $e_\L\DGp$ (provided that we have proved Conjecture \ref{conjmain}). Here $e_\L\in \DGp$ is a Heisenberg idempotent. Hence modulo Conjecture \ref{conjmain}, we have reduced the study of character sheaves in general to that in the Heisenberg case as above.

However, we have also proved that Conjecture \ref{conjmain} is equivalent to Conjecture \ref{conjweak}. Now Conjecture \ref{conjweak} in turn is just a statement about Heisenberg admissible pairs. Hence even to prove Conjecture \ref{conjmain} it is enough to consider the Heisenberg case.

Hence we are reduced to the study of the following case:\\
Let $(H,\L)$ be a Heisenberg admissible pair for an affine algebraic group $G$ and let $e_\L\in \DG$ be the corresponding Heisenberg idempotent.  We are then reduced to studying the braided triangulated category $e_\L\DG$. In particular to prove Conjecture \ref{conjmain}, we must show that $e_\L\DG$ has no non-trivial idempotents. Then we must define the notion of character sheaves in the category $e_\L\DG$ i.e. define the set $CS_{e_\L}(G)$. Then we may ask whether we can define this set purely in terms of the braided monoidal structure of $e_\L\DG$, namely does every braided autoequivalence of $e_\L\DG$ preserve the set $CS_{e_\L}(G)$?

Let $U$ be the unipotent radical of $G$. Then $(H,\L)$ is a Heisenberg admissible pair for the connected unipotent groups $U$ and by the results of \cite{De} $e\D_U(U)$ is equivalent to the bounded derived category pointed modular category corresponding to a metric group $(K_\L,\theta)$ (see \S\ref{anheisid}). Let $\Gamma=G/U$. Then $\Gamma^\circ$ is a connected reductive group. Using the observations from \S\ref{anheisid}, we may interpret the category $e_\L\D_U(G)$ as a twisted version of the category $\D(\Gamma)$ and the category $e_\L\DG$ as a twisted version of the category $\D_\Gamma(\Gamma)$.

\bex
Let $G$ be a reductive group and consider the Heisenberg admissible pair $(\{1\}, \Qlcl)$. (By Proposition \ref{admonred} this is the only admissible pair for $G$.) In this case the Heisenberg idempotent $e_{\{1\},\Qlcl}$ is the unit object $\delta_1\in \DG$. In this case the notion of character sheaves is defined in \cite{L}. However it is still an open question whether $\DG$ has any non-trivial idempotents. Also it is still an interesting question whether any braided autoequivalence of $\DG$ preserves the character sheaves on $G$.
\eex

\bex
Let us now consider the example of the so-called ``geometric Weil representation'' studied in \cite{GH}. Suppose that $\operatorname{char} \k\neq 2$. Let $V$ be a $\k$-vector space equipped with a symplectic form $B:V\times V\to \G_a.$ Let us form the Heisenberg group $U:=V\times \G_a$ with multiplication given by $(v,a)\cdot(w,b)=(v+w, a+b+\frac{1}{2}B(v,w))$. Let $\Gamma:=Sp(V)$, the symplectic group of $V$. Then $\Gamma$ acts on $U$ and we form the semidirect product $G=\Gamma\ltimes U$. Let $\L$ be any non-trivial multiplicative local system on $\G_a\subset U\subset G$. Then $(\G_a,\L)$ is a Heisenberg admissible pair for $G$. In this case, it is easy to check that $e\D_U(U)\cong D^b(\Vec)$. Moreover, using the results from \cite{GH} we can prove that in this case we have a braided equivalence $e_\L\DG\cong \D_\Gamma(\Gamma)$.
\eex

Once we analyze the Heisenberg case completely, we can use the results of this paper to define character sheaves in general.
\brk
Character sheaves on reductive groups as well as on unipotent groups are perverse sheaves (at least up to shift). We also expect character sheaves to be perverse in the {\it Heisenberg} case. However in general, character sheaves may not be perverse (even up to shift). This is because we will use induction functors to define character sheaves and these functors do not preserve perversity in general.
\erk

\appendix

\section{Appendix: Torus actions on connected unipotent groups}\label{torusactions}
\subsection{Equivariant multiplicative local systems}
Let us fix an embedding $\Q_p/\Z_p\hookrightarrow \Qlcl^*$. Now we can identify isomorphism classes of multiplicative $\Qlcl$-local systems on any connected (perfect) unipotent group $U$ and central extensions of $U$ by $\Q_p/\Z_p$. 
\blem\label{l:teqls}
Let $U$ be a connected unipotent group equipped with an action of a torus $T$ such that $U^T=\{1\}$.  Then there are no nontrivial $T$-equivariant multiplicative local systems on $U$.
\elem
\bpf
Let $\L$ be a $T$-equivariant multiplicative local system on $U$. It can be obtained from a central extension $0\to A\to \tilde{U}\to U\to 0$ with $\tilde{U}$ also connected and unipotent and $A$ a finite subgroup of $\Q_p/\Z_p$. Since central extensions of a connected group by a discrete group have no non-trivial automorphisms and since $\L$ is $T$-equivariant, we get an action of $T$ on $\tilde{U}$ making the above central extension $T$-equivariant. Taking $T$-fixed points, we get an exact sequence $0\to A\to \tilde{U}^T\to 1.$ Since $T$ is a torus and $\tilde{U}$ is a connected unipotent group, $\tilde{U}^T$ is connected (by \cite[\S18]{H}). Hence we must have $A\cong 0$, and hence $\L\cong \Qlcl.$
\epf

If a torus $T$ acts on a connected unipotent group $U$, then it follows from \cite[\S18]{H} that $U^T$ is connected. In this situation we have:
\blem\label{fp}
If $\L$ is a $T$-equivariant multiplicative local system on $U$ such that $\L|_{U^T}\cong \Qlcl$, then $\L\cong \Qlcl$. In other words, the restriction homomorphism $(U^*)^T\to (U^T)^*$ of perfect commutative unipotent groups is an injection.
\elem
\bpf
Let $\L$ be any $T$-equivariant multiplicative local system. As before, $\L$ comes from a $T$-equivariant central extension $0\to A\to \tU\to U\to 0$ with $\tU$ connected. Let us consider the central extension obtained by pullback $0\to A\to \t{U^T}\to U^T\to 0$.
 Let $\t{u}\in \t{U^T}$, that is $\t{u}$ maps to $U^T\subset U$. Hence $\t{u}\cdot t(\t{u}^{-1})\in A$ for every $t\in T$. Now since $T$ is connected and $A$ is a finite group, we must have $\t{u}\in \tU^T$. Hence we have $\t{U^T}=\tU^T$. Now let $\L$ be such that the restriction of $\L$ to $U^T$ is trivial. Hence the central extension $0\to A\to \t{U^T}=\tU^T\to U^T\to 0$ is trivial. On the other hand, $\tU^T$ is connected since $\tU$ is connected and $T$ is a torus. Hence we must have $A\cong 0$ or in other words $\L\cong \Qlcl.$
\epf

\blem
Let a reductive group $\Gamma$ act on a connected unipotent group $U$. Then we get the action of $\Gamma$ on the variety $U/{U^\Gamma}$. The only fixed point for this action is the trivial coset $U^\Gamma$.
\elem
\bpf
Suppose $u \in U$ is such that $c_\g(u)=u\cdot \g(u^{-1})\in U^\Gamma$ for each $\g\in \Gamma$. Since the action of $\Gamma$ on $U^\Gamma$ is trivial, we can check that $\g\mapsto c_\g(u)$ defines a group homomorphism $\Gamma\to U^\Gamma$ which must be trivial since $\Gamma$ is reductive and $U^\Gamma$ is unipotent. Hence $u\in U^\Gamma$. 
\epf


\subsection{Skewsymmetric biextensions equivariant under torus actions}

Let us continue to work in the world of perfect algebraic groups. Let $\cpu$ (resp. $\cpuc$) denote the category of perfect commutative unipotent groups (resp. perfect connected commutative unipotent groups). We refer to \cite[Appendix A.10]{B} for the theory of skewsymmetric isogenies and biextensions of perfect connected commutative unipotent groups and their associated metric groups. Let $V\in\cpuc$ and let $V\stackrel{\phi}{\rightarrow}V^*$ be a skewsymmetric isogeny with kernel $K_\phi$. According to \cite{Da},\cite{B}, there is the metric group $(K_\phi,\theta)$ associated with the skewsymmetric isogeny $\phi.$ Let a torus $T$ act on $V$ and hence on $V^*$. Suppose that the isogeny $\phi$ is $T$-equivariant, i.e. we have a $T$-equivariant exact sequence 
\beq\label{phi}
0\to K_\phi\to V\stackrel{\phi}{\rightarrow}V^*\to 0.
\eeq

Since $T$ is connected, the action of $T$ on $K_\phi$ is trivial and hence $K_\phi\subset V^T$. Again since $T$ is connected and $K_\phi$ is finite, taking $T$-fixed points of (\ref{phi}) (and using the argument used in proof of Lemma \ref{fp}), we get the exact sequence 
\beq\label{fpsksymm}
0\to K_\phi\to V^T\to (V^*)^T\to 0.
\eeq
Note that according to Lemma \ref{fp}, the restriction of multiplicative local systems induces an injection $(V^*)^T\hookrightarrow (V^T)^*.$ Moreover, by the above isogeny, we have the equality $\dim((V^*)^T)=\dim(V^T)=\dim((V^T)^*)$. Hence the above inclusion must in fact be an isomorphism $(V^*)^T\stackrel{\cong}{\rightarrow} (V^T)^*$. In other words, the quotient $(V^T)^*$ of $V^*$ can be identified with the subgroup $(V^*)^T\subset V^*$. Hence the short exact sequence $0\to (V/V^T)^*\to V^*\to (V^T)^*\cong (V^*)^T\to 0$ splits as $V^*\cong (V^*)^T\oplus (V/V^T)^*.$ Taking duals, we see that the short exact sequence $0\to V^T\to V\to V/V^T\to 0$ also splits as $V\cong V^T\oplus V/V^T.$ We see from (\ref{fpsksymm}) (and the identification $(V^*)^T\cong (V^T)^*$) that the restriction of $\phi$ to $V^T$ is again a skewsymmetric isogeny which has the same metric group $(K_\phi,\theta)$ associated with it. In fact, we see that we have the following:
\bprop\label{skewsymm}
Let $V\in \cpuc$ be equipped with the action of a torus $T$. Let $V\stackrel{\phi}{\rightarrow} V^*$ be a $T$-equivariant skewsymmetric isogeny with associated metric group $(K_\phi,\theta).$ Then $V$ splits as $V^T\oplus V'$ where $V'\subset V$ is $T$-invariant and the isogeny $\phi$ splits into:
\bit
\item a skewsymmetric isogeny $\phi^T:V^T\to (V^T)^*$ with associated metric group $(K_\phi,\theta)$ and
\item a $T$-equivariant skewsymmetric \emph{isomorphism} $\phi':V'\to (V')^*.$
\eit
\eprop
In this decomposition, $T$ acts trivially on $V^T$ and $V'$ has no nonzero $T$-fixed points.

\subsection{Torus actions on vector groups}
\bdefn\label{defvecgp}
A group $V\in \cpuc$ is said to be a vector group if there is an isomorphism $V\xto{\cong} \G_a^n$ for some non-negative integer $n.$ We can use such an isomorphism to define an action of $\G_m$ on $V$ using its natural action on $\G_a^n$. This action is defined to be {\it a linear structure} on the vector group $V$. If $V,W$ are vector groups equipped with a linear structure, we say that a group homomorphism $V\to W$ is linear if it commutes with the chosen $\G_m$-scalings on $V,W$.
\edefn

\brk
Equivalently, $V\in \cpuc$ is a vector group if and only if $p\cdot V=0$.
\erk

The following result is well known (cf. \cite[\S4]{C})
\bthm\label{linearize}
Let $T$ be a torus acting on a vector group $V$ by group automorphisms. Then there exists a linear structure on $V$ such that $T$ in fact acts by linear automorphisms with respect to this linear structure.
\ethm

Hence in order to study vector groups equipped with torus actions, it is enough to study linear actions. If $V$ is such a linear representation of $T$, then we have a decomposition of $V$ into weight spaces $V=\bigoplus{V_\lambda}$. If $V^*$ is the contragradient representation, then we have 
$$V^*=\bigoplus{(V_\lambda)^*}=\bigoplus{(V^*)_{-\l}}$$ gives the decomposition of $V^*$ into weight spaces since we have an identification $(V^*)_{-\l}=(V_\l)^*$. However, if we have two linear representations $V,W$ of $T$, we must also study {\it non-linear} $T$-maps  $V\xto{\phi} W$. It is clear that $\phi(V_0)=\phi(V^T)\subset W_0=W^T$. We still need to see what happens to the other weight spaces.

For a weight $0\neq\lambda$ of $V$ and $\lambda'$ of $W$, let us consider the induced $T$-map $V_\lambda\xto{\phi'} W_{\lambda'}$. For $v\in V_\l, t\in T$ we have
\beq\label{tact}
\phi'(\l(t)v)=\phi'(tv)=t\phi'(v)=\l'(t)\phi'(v).
\eeq
Since $\l\neq 0$, this implies that the line spanned by $v$ maps to the line spanned by $\phi'(v)$. Hence if $\phi'$ is non-zero, then $\Hom_T(\G_{a,\l}, \G_{a,\l'})\neq 0,$ where $\G_{a,\l}, \G_{a,\l'}$ are the one dimensional representations of $T$ corresponding to the weights $\l,\l'$ respectively.

Let $\phi':\G_{a,\l}\to \G_{a,\l'}$ be a non-zero $T$-map with $\l\neq 0$. We can rescale this map so that $\phi'(1)=1.$ Under the inclusion $\G_m\subset \G_a$, we see that $\phi'(\l(t))=\l'(t)$ by putting $v=1$ in (\ref{tact}). Hence we must have $\phi'=\tau^n$ (where $\tau:\G_a\to \G_a$ is the Frobenius automorphism) for some (possibly negative) $n\in \Z$ and $\l'(t)=\l(t)^{p^n}$, i.e. $\l'=p^n\l$ in usual additive notation. 

Hence we see that $\Hom_T(\G_{a,\l}, \G_{a,p^n\l})= \{c\tau^n|c\in \k\}$ if $\l\neq0$,  $\Hom_T(\G_{a,0}, \G_{a,0})=\Hom(\G_{a}, \G_{a})=k\{\tau,\tau^{-1}\}$ and that all other $\Hom$'s are zero.

For a weight $\l:T\to \G_m$, and a linear representation $V$ of $T$, let $V_{\l}^{(p)}=\sum\limits_{n\in \Z}{V_{p^n\l}}\subset V$. We have $V_0^{(p)}=V_0$. We have proved
\bthm\label{Tmap}
Let $V$ be a linear representation of a torus $T$. Then we have a decomposition $$V=\bigoplus{V_\l^{(p)}}$$ into $T$-stable subspaces. For the contragradient representation the decomposition is given by 
$$V^*=\bigoplus{(V^*)_{-\l}^{(p)}}=\bigoplus{\left(V_{\l}^{(p)}\right)^*}.$$ 
If $V,W$ are linear representations of $T$ and if $\phi:V\to W$ is a (not necessarily linear) $T$-map, then $\phi(V_\l^{(p)})\subset W_\l^{(p)}$. 
\ethm

\subsection{$T$-stable maximal isotropic subgroups}\label{equivisotropic}
Let $V\in \cpuc$ and let $\phi:V\to V^*$ be a skew-symmetric biextension (not necessarily an isogeny). For a closed connected subgroup $L\subset V$, let $L^\perp\subset V$ be the kernel of the composition $V\xto{\phi}V^*\to L^*$ and $L^{\perp\circ}$ its neutral connected component. For example $V^\perp=\ker(\phi).$ The biextension $\phi$ is an isogeny if and only if $V^{\perp\circ}=0$.

\bdefn\label{defisotropic}
We say that a connected subgroup $L\subset V$ is isotropic if $L\subset L^\perp$, or equivalently if $L\subset L^{\perp\circ}$. If $L$ is an isotropic subgroup, then $\phi$ induces a skew-symmetric biextension on the subquotient
\beq
\phi_L: L^{\perp\circ}/L \to  (L^{\perp\circ}/L)^*.
\eeq 
A maximal isotropic subgroup is an isotropic subgroup that is maximal among all isotropic subgroups. If $L$ is maximal isotropic, then we must have $\dim(L^{\perp\circ}/L)\leq 1$. (See \cite[Appendix A.2]{De13}.)
\edefn

Now suppose that $V$ is equipped with an action of a torus $T$  such that the biextension $\phi:V\to V^*$ is $T$-equivariant. Our goal is to prove
\bthm\label{stabisot}
In the situation above, there exists a maximal isotropic subgroup $L\subset V$ such that $L$ is $T$-stable.
\ethm
We carried out a similar argument in \cite[Appendix A.2]{De13}. As in {\it loc. cit.}, if we have a $T$-stable isotropic subgroup we can pass to the associated subquotient. Moreover, the assertion is obvious in case $\dim(V)\leq 1.$ Hence it only remains to prove:
\bprop\label{dimgeq2}
In the above notation, if $\dim(V)\geq 2$, then there exists a non-trivial $T$-stable isotropic subgroup.
\eprop
\bpf
By \cite[Lem. A.9]{De13}, we easily reduce to the case $p\cdot V=0$, i.e. $V$  a vector group. Let us choose a linear structure on $V$ such that the action of $T$ is linear. Hence the action of $T$ on $V^*$ is also linear and we can think of $V^*$ as the contragradient linear representation. By Theorem \ref{Tmap}, we have  decompositions $$V=\bigoplus{V_\l^{(p)}},$$
$$V^*=\bigoplus{\left(V_\l^{(p)}\right)^*=\bigoplus{(V^*)_{-\l}^{(p)}}}$$ and inclusions $\phi\left(V_\l^{(p)}\right)\subset \left(V_{-\l}^{(p)}\right)^*$. Hence if $\l\neq 0,$ then $V_{-\l}^{(p)}$ is a ($T$-stable) isotropic subgroup. Hence if $V$ has non-zero weights, then we get a non-trivial $T$-stable isotropic subgroup. On the other hand if $T$ acts trivially on $V$, then any isotropic subgroup is $T$-stable. Since $\dim(V)\geq 2$ we know that $V$ must have a non-trivial isotropic subgroup (see \cite[Appendix A.10]{B} or \cite[Prop. A.10]{De13}). 

\epf

\bibliographystyle{ams-alpha}

\end{document}